\pdfoutput=1
\documentclass[11pt]{article}
\usepackage[margin=1in]{geometry}
\usepackage[T1]{fontenc}
\usepackage{amsmath,amssymb,amsthm}
\usepackage{graphicx}
\usepackage{booktabs}
\usepackage{enumitem}
\usepackage{multirow}
\usepackage[table]{xcolor}
\usepackage{subcaption}
\usepackage{placeins}
\usepackage[raggedright]{titlesec}
\usepackage[title]{appendix}
\usepackage[round,authoryear]{natbib}
\usepackage[hidelinks]{hyperref}

\DeclareMathOperator*{\Minimize}{Minimize}
\newlength{\fullwidth}
\AtBeginDocument{\setlength{\fullwidth}{\textwidth}}

\newcommand{\TABLE}[3]{\centering #2\par\caption{#1 #3}}
\newcommand{\FIGURE}[3]{\centering #1\caption{#2 #3}}

\title{Using block coordinate descent for solving\\ the close enough traveling salesman problem\thanks{This work has been partially supported by the Brazilian agencies FAPESP (grants 2022/05803-3 and 2023/08706-1) and CNPq (grant 302073/2022-1).}}

\author{
    Ernesto G. Birgin\thanks{Department of Computer Science, Institute of Mathematics, Statistics and Computer Science, University of S\~ao Paulo, Rua do
    Mat\~ao, 1010, Cidade Universit\'aria, 05508-090, S\~ao Paulo, SP,
    Brazil. e-mail: egbirgin@ime.usp.br}
\and 
    Paula C. R. Ertel\thanks{Department of Applied Mathematics, Institute of Mathematics, Statistics and Computer Science, University of S\~ao Paulo, Rua do
    Mat\~ao, 1010, Cidade Universit\'aria, 05508-090, S\~ao Paulo, SP,
    Brazil. e-mail: \{paulaertel|rogerio.matos\}@ime.usp.br}
\and 
    Rogério S. Matos\footnotemark[3]
}

\date{November 24, 2026}

\begin{document}

\maketitle

\begin{abstract}
The Close Enough Traveling Salesman Problem (CETSP) is a variant of the Euclidean Traveling Salesman Problem in which the salesperson must visit some point within a disk associated with each customer. Evaluating a visiting sequence requires determining an optimal visit point within each disk. This building block of virtually every CETSP algorithm is typically computed either by solving a second-order cone programming (SOCP) problem with a general-purpose solver, which can be expensive, or by discretizing the disks, which sacrifices optimality. We propose a block coordinate descent (BCD) method that computes these visit points exactly, without discretization, optimizing one visit point at a time. We show that the resulting single-disk subproblem reduces to a one-dimensional minimization, and we characterize its stationary points through the roots of a quartic polynomial. Of the three resulting solution strategies, which agree up to machine precision, Newton's method is the most efficient. Compared with an off-the-shelf SOCP solver, BCD evaluates fixed tours more than an order of magnitude faster. Embedded in an iterated local search (ILS-BCD), it finds the known optimal solutions of 720 small benchmark instances and, on 62 large ones, outperforms all heuristics of their original study in 47 and matches the best of them in 11 others. Within the memetic algorithm MA-CETSP, replacing the general-purpose SOCP solver with BCD and adding an exact local search based on the single-disk solver yields a method that, over the best of 20 runs, is better than the original algorithm in 20 of the 62 large instances and equal in the others, and improves on its best published values in 21. BCD is a fast and exact building block for evaluating visiting sequences that can replace general-purpose SOCP solvers within CETSP algorithms and enables searches that would otherwise be computationally prohibitive.

\medskip\noindent
\textbf{Keywords:} Close Enough Traveling Salesman Problem, block coordinate descent, single-disk problem, iterated local search, memetic algorithm.
\end{abstract}

\section{Introduction}

The Close Enough Traveling Salesman Problem (CETSP) is a variant of the Euclidean Traveling Salesman Problem in which each customer is associated with a disk rather than a point. The objective is to determine both a visiting sequence and a visit point within each disk so as to minimize the total length of the resulting tour. The problem that is now known as the CETSP was introduced, without being given this specific name, in~\citet{Arkin1994}. The problems considered in that paper are referred to as Geometric Covering Salesman Problems. The authors consider settings in which the neighborhoods to be visited are parallel line segments, or translations of a convex neighborhood, such as a disk or a rectangle. They propose approximation algorithms for solving these problems. Thus, to the best of our knowledge, the CETSP, as a problem, was first introduced in~\citet{Arkin1994}, although the term CETSP was coined only later. A related, but different, problem was considered in~\citet{Current1989}, under the name Covering Salesman Problem (CSP). In the CSP, a set of nodes is given, and the objective is to find a shortest tour visiting a subset of these nodes such that every node not included in the tour lies within a prescribed distance of at least one node on the tour. A similar problem, known as the Covering Tour Problem (CTP), was later considered in~\citet{Gendreau1997}. In the CTP, the vertices of a graph are divided into two sets: optional vertices, which may be visited, and required vertices, which must be covered, meaning that they must lie within a prescribed distance of a visited vertex. Assuming that edge weights satisfy the triangle inequality, the objective is to find a minimum-length Hamiltonian cycle over a subset of the optional vertices that covers all required vertices. Recent literature also addresses a related variant known as the Carrier-Vehicle TSP~\citep{erdogan2021,Vieira2026}. In this framework, a slow carrier (e.g., a vessel) and a fast vehicle (e.g., a helicopter or a drone) execute emergency response, post-disaster relief, aerial inspection, environmental mapping, continuous tracking, or payload transport operations. The carrier needs only to come close enough to each target, enabling the fast vehicle to launch, visit the target, and return to the moving base within its operating range. Consequently, optimal launch and recovery points must be determined in continuous space to minimize travel time or distance for either or both vehicles. A traveling salesman problem with specific neighborhoods was also considered in~\citet{Dumitrescu2003}, where approximation algorithms were developed for the case of identical disks. The term Close Enough Traveling Salesman Problem (CETSP) was coined in~\citet{Gulczynski2006}, where the problem was studied in the form in which it is commonly considered today. In that work, the problem arises in the context of radio-frequency identification (RFID) technology, which allows utility companies to read meters from a distance. Identical disks are considered as neighborhoods.

The CETSP has attracted considerable attention because it arises naturally in applications where a target needs to be approached only within a prescribed distance. Such applications include surveillance, wireless communication, data collection, and unmanned aerial vehicle routing. For example, in~\citet{Yuan2007}, distributed sensors are sparsely deployed in the Euclidean plane, and a mobile robot must reach each of them within a disk centered at the sensor, whose radius is the maximum distance at which its data can be collected. A related application is considered in~\citet{Shuttleworth2008}, where radio-frequency identification enables automated meter reading at a distance, so that a meter reader need only pass within a given radius of each customer. There, the route is constrained to a realistic street network. In~\citet{Poikonen2017}, a vehicle routing problem involving multiple trucks and drones is considered. In this last-mile delivery setting, the trucks need only to approach the delivery locations within a certain distance, at which point the drones are launched to complete the deliveries. The allowable distance between a truck's visit point and a delivery location is related to the drones' battery life and, consequently, to their operational range. Finally, \citet{diplacido2022cetsp} models as a CETSP a real-world application in which a drone equipped with a thermal camera performs diagnostic reconnaissance tours over a photovoltaic field, with the aim of assessing the operating condition of the solar panels.

Solving the CETSP involves two interrelated decisions: determining a visiting sequence and, for a given sequence, determining a visit point within each disk. The methods proposed for solving the CETSP can be classified according to how they determine the actual visit point within each disk. A first family of approaches discretizes each neighborhood into a finite set of candidate points and formulates the resulting problem as a mixed-integer program or performs a heuristic search over this finite set. In~\citet{Behdani2014}, each disk is discretized into five boundary points corresponding to the endpoints of four consecutive arcs. These points are then used to construct a solution, while the arcs are used to derive a lower bound on the optimal tour length. This scheme is refined in~\citet{Carrabs2017b,Carrabs2020}, where an adaptive discretization with up to eight points per disk is proposed, with the points being placed at the midpoints of the sides of a regular polygon contained in the disk rather than on its circumference. A similar discretization, in which four equally spaced points on the boundary of each disk are considered, is used in~\citet{Nedjati2020}. There, the discretization is used to construct a visiting sequence, after which the visit points are computed using an off-the-shelf global optimization solver. In~\citet{Dong2007}, the problem is first formulated as a mixed-integer nonlinear program, after which heuristics that combine discretization of the set of candidate visit points with clustering strategies are proposed to solve it in practice. In~\citet{Yuan2007}, a finite set of angles is considered, corresponding to a discretization of the disk boundary. More recently,~\citet{Fan2025} proposed a deep reinforcement learning approach that constructs the visiting sequence incrementally. Starting from an empty sequence, the method simultaneously selects the next customer to be visited and its visit point at each step, with the latter determined through the discretization scheme proposed in~\citet{Carrabs2017b}.

A second family of approaches instead solves, for a fixed visiting sequence, a continuous Second-Order Cone Program (SOCP) to determine the visit points exactly. This idea originates in~\citet{Mennell2011}, where the solution procedure consists of three steps. First, specific points within so-called Steiner zones --- regions formed by the intersection of one or more disks --- are identified. These points are then used to define a visiting sequence. Finally, given this sequence, the visit points are optimally repositioned by solving an SOCP. In~\citet{Coutinho2016}, such an SOCP subproblem is embedded within a branch-and-bound search over visiting sequences, yielding an exact algorithm for the CETSP; this branch-and-bound framework is later improved in~\citet{Zhang2023}. An SOCP is also used as a post-optimization step, applied to a visiting sequence obtained in different ways. In~\citet{Carrabs2017} and~\citet{Carrabs2020}, a discretization-based model is used to produce a solution consisting of a visiting sequence and the corresponding visit points, and an SOCP is then used to optimally recompute the visit points for the given sequence, thus potentially improving the solution. In~\citet{diplacido2022cetsp}, an SOCP subproblem is likewise solved for a randomly selected individual in each generation of a genetic algorithm, rather than for a single fixed visiting sequence. In~\citet{erdogan2021}, only a limited number of SOCP subproblems are solved within an Iterated Local Search (ILS) heuristic framework, as the authors explicitly note that repeatedly solving SOCPs inside local search operators would be computationally prohibitive. A related strategy is adopted in~\citet{Vieira2026}, where the conic formulation is replaced by a linearization whose approximation error is corrected by an exact refinement procedure.

While discretization-based methods sacrifice optimality by restricting the visit points to a finite set of candidates, methods based on solving an SOCP determine the optimal visit points for a given visiting sequence. However, solving an SOCP for every fixed sequence can be computationally expensive. This has motivated more recent heuristics that exploit the nonlinear structure of the fixed-sequence problem to determine the visit points without solving the full SOCP. However, these methods treat the determination of each visit point individually as a continuous optimization problem and solve each such problem approximately. Thus, they still sacrifice optimality in two ways: by not treating the determination of the visit points as an integrated optimization problem and by solving each individual problem only approximately. In~\citet{Wang2019}, Steiner zones are used, and visit points within the zones are computed by an iterative scheme that alternates between updating the even- and odd-indexed points until no further improvement is obtained. With the preceding and succeeding visit points fixed, the problem of positioning a visit point is reduced to a univariate nonlinear optimization problem, which is approximately solved using bisection to find a zero of its derivative. Similarly,~\citet{Lei2024} recognizes the underlying nonlinear structure of the problem and avoids the computational cost of solving an SOCP while exploring visiting sequences. Within its local search, the visit points are computed heuristically, and an SOCP is solved only for the sequence returned at the end of each local search. In~\citet{Yang2018}, a hybrid heuristic method is proposed in which particle swarm optimization searches the visit points in an outer loop, while a genetic algorithm determines the visiting sequence for each candidate set of points in an inner loop.

Most CETSP algorithms repeatedly generate candidate visiting sequences. Because the total tour length depends on the specific visit points, evaluating each candidate sequence requires solving a continuous optimization problem to determine the optimal visit points within the disks. This fixed-sequence problem is therefore a fundamental computational primitive, as it must be solved whenever a candidate sequence is evaluated. The problem is convex and can be formulated as a second-order cone program (SOCP). As highlighted above, most works in the literature either formulate and solve this problem directly as an SOCP, typically using an off-the-shelf general-purpose optimization solver, or replace it with an approximate problem based on a coarse discretization and solve the resulting problem exactly or heuristically.

Block coordinate descent (BCD) is one of the most natural approaches for exploiting the structure of this evaluation problem. At each iteration, one block of variables is optimized while all remaining variables are kept fixed. In the present context, each block consists of the two coordinates of the visit point associated with one disk. Therefore, every BCD iteration reduces to determining the optimal location of a single visit point while the remaining visit points are fixed. Furthermore, due to the nature of the objective function, which corresponds to the sum of distances between consecutive tour points, the optimal placement of each visit point depends solely on its immediate predecessor and successor.
BCD was used to evaluate visiting sequences in two earlier works, both limited in different ways. In~\citet{Birgin2022}, the neighborhoods are not disks but polygons and each block update solves a regularized subproblem for which no efficient exact method is known; it is solved approximately, by discretization and brute force. In~\citet{Ertel2023}, which considers non-overlapping disks, the problem of locating a single visit point is solved by a globalized Newton method. Globalization techniques are known to slow down methods with fast local convergence, such as Newton's method, since they fall back on slower steps whenever the Newton step is not acceptable~\citep[\S2.5 and \S6.2]{DennisSchnabel1996}, and the method therefore becomes inefficient. Furthermore, globalized Newton methods only guarantee convergence to stationary points, which may not be global minimizers. In both works, the combinatorial part is straightforward: a constructive heuristic followed by a single local search, with every neighbor evaluated by a complete run of BCD, and neither method was compared with others on standard benchmarks. The present work removes these limitations. Its main contribution is the development of a specialized BCD algorithm for solving the fixed-sequence evaluation problem arising in the Euclidean CETSP. By exploiting the geometry of the one-disk subproblem, we show that each BCD update can be reduced to the minimization of a nonlinear function of a single variable. This reduction leads to a highly efficient implementation based on Newton's method, eliminating the need for a general-purpose optimization solver.

To provide an independent and robust assessment of the proposed approach, we also derive an algebraic characterization of all stationary points of the one-disk subproblem. Starting from the first-order optimality conditions, we derive a quartic polynomial whose real roots are in one-to-one correspondence with the stationary points of the subproblem. This characterization leads to two alternative numerical procedures based on polynomial root finding: the Durand--Kerner method and a companion-matrix eigenvalue formulation. These procedures provide independent means of computing the stationary points and allow us to assess the proposed one-dimensional optimization approach both in terms of solution quality and computational efficiency. Numerical experiments show that all three approaches produce the same solutions up to machine precision, providing a robust validation of the proposed method. Newton's method is preferable though, because it consistently requires substantially less computational effort than the quartic-based approaches.

When compared to a conventional off-the-shelf solver, the BCD approach evaluates visiting sequences more than an order of magnitude faster and finds the same solutions but at a significantly higher precision. Therefore BCD establishes a fast and exact building block for CETSP algorithms. To showcase its effectiveness, we integrate BCD into an ILS (ILS-BCD). ILS-BCD identified all global minimizers across 720 small-scale instances and outperformed the heuristics evaluated in~\citet{mennell2009thesis} on 62 large-scale instances, obtaining a better solution than all of them in 47 instances and matching the best of them in 11 others. We also integrate BCD into the memetic algorithm MA-CETSP~\citep{Lei2024}, one of the best-performing methods for these instances, replacing its general-purpose SOCP solver and using the single-disk solver within its local search. On the 62 large instances, the resulting method, MA-CETSP-BCD, is better than MA-CETSP in 20 instances and equal in the other 42 when the best of 20 runs of each method is compared, and better than the best values reported in~\citet{Lei2024}, which were obtained in separate runs, in 21 instances and equal in the other 41. An ablation study shows that replacing the solver alone preserves the quality of MA-CETSP while reducing its running time, and that the improvements come from the uses of BCD and of the single-disk solver that its speed and accuracy make possible.

The remainder of the paper is organized as follows. Section~\ref{sec2} describes the CETSP problem and presents its mathematical formulation. Section~\ref{sec3} introduces the proposed BCD method, the reduction of the single-disk subproblem to a one-dimensional optimization problem, and its algebraic characterization via a quartic polynomial. Section~\ref{sec4} reports numerical experiments comparing the alternative subproblem solvers and evaluating BCD against CPLEX. Section~\ref{secmeta} assesses BCD within two metaheuristics: the trajectory-based ILS-BCD and the population-based MA-CETSP-BCD. Finally, Section~\ref{sec5} concludes the paper.

\section{Problem formulation} \label{sec2}

In this section, we describe the problem through its mathematical mixed integer quadratically constrained programming (MIQCP) formulation, which provides an unambiguous description of the problem under consideration. We then consider a simplified version of this formulation, in which the variables determining the visiting sequence are fixed. This reduced model is a second-order cone programming (SOCP) problem whose solution gives the optimal visit points for a given visiting sequence, which is precisely the problem addressed by the proposed BCD method. It allows us to use a conventional off-the-shelf solver to compare results, and thus evaluate the correctness and effectiveness of the proposed method.

We consider a set of $n$ customers to be visited, where each customer $i \in \{1, \dots, n\}$ is characterized by a circular neighborhood defined by its fixed center coordinates $c_i = (c_i^x, c_i^y) \in \mathbb{R}^2$ and an associated coverage radius $r_i \ge 0$. The presence of a depot can be incorporated into the model by treating it as a customer whose circular neighborhood has zero radius. To guarantee the proper relaxation of the distance constraints, we also define a sufficiently large constant $M$ (Big-$M$), set to be at least equal to the maximum possible distance between any two points within the customer neighborhoods.

The decision variables of the model are the spatial coordinates $y_i = (y_i^x, y_i^y) \in \mathbb{R}^2$ of the visit point of customer $i \in \{1, \dots, n\}$; the binary variables $x_{ij} \in \{0, 1\}$, equal to 1 if the tour goes directly from customer $i$'s neighborhood to customer $j$'s neighborhood, and 0 otherwise ($i, j \in \{1, \dots, n\}, i \neq j$); the continuous variables $\xi_{ij} \ge 0$, representing the Euclidean distance between $y_i$ and $y_j$ ($i, j \in \{1, \dots, n\}, i \neq j$); the continuous variables $\gamma_{ij} \ge 0$, representing the cost of traversing directly from customer $i$'s neighborhood to customer $j$'s neighborhood, counted in the objective only when $x_{ij} = 1$; and the continuous auxiliary variables $u_i \in \mathbb{R}$ for subtour elimination ($i \in \{1, \dots, n\}$). The model follows.

\begin{align}
\Minimize & \sum_{i=1}^n \sum_{\substack{j=1 \\ j \neq i}}^n \gamma_{ij} \label{eq:obj} \\
\text{s.t.} 
& \sum_{\substack{j=1 \\ j \neq i}}^n x_{ij} = 1 & i \in \{1, \dots, n\} \label{eq:degree_out} \\
& \sum_{\substack{i=1 \\ i \neq j}}^n x_{ij} = 1 & j \in \{1, \dots, n\} \label{eq:degree_in} \\
& u_i - u_j + n x_{ij} \le n - 1 & i, j \in \{2, \dots, n\}, i \neq j \label{eq:mtz} \\
& (y_i^x - c_i^x)^2 + (y_i^y - c_i^y)^2 \le r_i^2 & i \in \{1, \dots, n\} \label{eq:neighborhood} \\
& \xi_{ij}^2 \ge (y_i^x - y_j^x)^2 + (y_i^y - y_j^y)^2 & i, j \in \{1, \dots, n\}, i \neq j \label{eq:distance} \\
& \gamma_{ij} \le M x_{ij} & i, j \in \{1, \dots, n\}, i \neq j \label{eq:gamma_upper} \\
& \gamma_{ij} \ge \xi_{ij} - M(1 - x_{ij}) & i, j \in \{1, \dots, n\}, i \neq j \label{eq:gamma_lower} \\
& x_{ij} \in \{0, 1\} & i, j \in \{1, \dots, n\}, i \neq j \label{eq:var_x} \\
& \xi_{ij} \ge 0 & i, j \in \{1, \dots, n\}, i \neq j \label{eq:var_xi} \\
& \gamma_{ij} \ge 0 & i, j \in \{1, \dots, n\}, i \neq j \label{eq:var_gamma} \\
& 1 \le u_i \le n & i \in \{1, \dots, n\} \label{eq:var_u}
\end{align}

The objective function \eqref{eq:obj} minimizes the total tour length by taking the sum of the arc costs $\gamma_{ij}$ over all customer pairs. Constraints \eqref{eq:degree_out} and \eqref{eq:degree_in} guarantee tour continuity by requiring that the tour enters and leaves the neighborhood of each customer $i \in \{1, \dots, n\}$ exactly once, while constraints \eqref{eq:mtz} represent the Miller--Tucker--Zemlin (MTZ) formulation that prevents subtours, using customer $1$ as an arbitrary reference node since there is no depot. Spatial feasibility is maintained by constraints \eqref{eq:neighborhood}, which enforce that the visit point coordinates $(y_i^x, y_i^y)$ chosen for customer $i$ lie within the circular region bounded by center $(c_i^x, c_i^y)$ and radius $r_i$. Constraints \eqref{eq:distance} are second-order cone constraints defining $\xi_{ij}$ as an upper bound on the Euclidean distance between $y_i$ and $y_j$; since the objective only rewards keeping the $\gamma_{ij}$, and through them the $\xi_{ij}$, as small as possible, this bound is tight at any optimal solution. The link between $\xi_{ij}$ and the arc cost $\gamma_{ij}$ actually counted in the objective is established by constraints \eqref{eq:gamma_upper}--\eqref{eq:gamma_lower} via a Big-$M$ formulation: when $x_{ij} = 0$, constraint \eqref{eq:gamma_upper} forces $\gamma_{ij} = 0$, while constraint \eqref{eq:gamma_lower} becomes vacuous; when $x_{ij} = 1$, constraint \eqref{eq:gamma_upper} becomes vacuous, while constraint \eqref{eq:gamma_lower} forces $\gamma_{ij} \ge \xi_{ij}$, which the minimization objective drives to equality. Finally, constraints \eqref{eq:var_x}--\eqref{eq:var_u} set the domains and bounds for the arc, distance, cost, and sequence variables.

Fixing a visiting sequence defined by a permutation $\pi$ of $\{1,\dots,n\}$ amounts to setting $x_{\pi(i-1),\pi(i)}=1$ for $i=1,\dots,n$, with $\pi(0)=\pi(n)$ and all other arc variables to zero. (For later use, we also define $\pi(n+1)=\pi(1)$.) Once the tour is fixed, the arc and subtour-elimination constraints~\eqref{eq:degree_out}--\eqref{eq:mtz} can be omitted. Moreover, only the distances associated with the selected arcs need to be represented, making the variables $\gamma_{ij}$ and the Big-$M$ constraints~\eqref{eq:gamma_upper} and~\eqref{eq:gamma_lower} unnecessary. Defining $t_i=\xi_{\pi(i-1),\pi(i)}$ for $i=1,\dots,n$ yields the second-order cone programming (SOCP) problem
\begin{equation} \label{socp-fixed-sequence}
\begin{array}{rll}
\displaystyle \Minimize_{\substack{y_1,\dots,y_n \in \mathbb{R}^2 \\[0.7mm] t \in \mathbb{R}^n}} & \displaystyle\sum_{i=1}^{n} t_i & \\[6mm]
\mbox{subject to} & \|y_{\pi(i)}-y_{\pi(i-1)}\|_2 \leq t_i, & i=1,\dots,n, \\[1mm]
& \|y_i-c_i\|_2 \leq r_i, & i=1,\dots,n.
\end{array}
\end{equation}
At any optimal solution, each $t_i$ equals the corresponding Euclidean distance. Thus, this SOCP formulation is equivalent to minimizing the total tour length directly over the visit points for the fixed visiting sequence.

\section{Block coordinate descent for evaluating a visiting sequence} \label{sec3}

Block Coordinate Descent (BCD) is an iterative optimization framework that partitions the decision variables into smaller subsets, or blocks. At each iteration, the method minimizes the objective function with respect to one block while holding all other blocks fixed. This procedure reduces a high-dimensional optimization problem to a sequence of lower-dimensional subproblems that can often be solved efficiently or even in closed form. Blocks may be selected cyclically, randomly, or according to a greedy rule, such as a Gauss--Southwell rule. BCD is particularly attractive when the constraints are separable across blocks and the individual block subproblems are substantially easier to solve than the full problem.

In this section, we develop a specialized BCD method for computing the visit points for a fixed visiting sequence. Eliminating the auxiliary variables $t\in\mathbb{R}^n$ from problem~\eqref{socp-fixed-sequence} yields the equivalent formulation
\begin{equation} \label{probfixedseq}
\begin{array}{rl}
\displaystyle\Minimize_{y_1,\dots,y_n \in \mathbb{R}^2} & \displaystyle\sum_{i=1}^{n} \|y_{\pi(i)}-y_{\pi(i-1)}\|_2 \\[5mm]
\mbox{subject to} & \|y_i-c_i\|_2 \leq r_i, \quad i=1,\dots,n.
\end{array}
\end{equation}
BCD sequentially updates each visit point $y_{\pi(i)}$ while holding all other visit points fixed. Since only the two incident tour edges depend on $y_{\pi(i)}$, each block update reduces to finding a point in the disk $B(c_{\pi(i)},r_{\pi(i)})=\{u\in\mathbb{R}^2:\|u-c_{\pi(i)}\|_2\leq r_{\pi(i)}\}$ that minimizes the sum of its Euclidean distances to the fixed neighboring visit points $y_{\pi(i-1)}$ and $y_{\pi(i+1)}$. The resulting single-disk subproblem is the two-dimensional convex optimization problem
\begin{equation} \label{singlediskprob}
\begin{array}{rl}
\displaystyle\Minimize_{y_{\pi(i)}\in\mathbb{R}^2} & \|y_{\pi(i)}-y_{\pi(i-1)}\|_2 + \|y_{\pi(i)}-y_{\pi(i+1)}\|_2 \\[4mm]
\mbox{subject to} & y_{\pi(i)} \in B(c_{\pi(i)},r_{\pi(i)}).
\end{array}
\end{equation}

The proposed BCD sweeps cyclically through the sequence of neighborhoods, updating each visit point $y_{\pi(i)}$ ($i=1,\dots,n$) individually while keeping all other coordinates fixed. The iterative process starts from given locations $y_{\pi(i)}^{(0)} \in B(c_{\pi(i)},r_{\pi(i)})$ for $i=1,\dots,n$. At the $i$-th step of cycle~$k$, it uses the current visit points $y_{\pi(i-1)}^{(k)}$ and $y_{\pi(i+1)}^{(k-1)}$ of the adjacent neighbors to obtain a new visit point $y_{\pi(i)}^{(k)} \in B(c_{\pi(i)},r_{\pi(i)})$ for customer $\pi(i)$. This cyclic process repeats until the relative decrease of the tour length in a full cycle drops below the tolerance $\epsilon_{\textsc{bcd}}=\sqrt{\epsilon_{\text{mach}}}\approx 1.49\times10^{-8}$, where $\epsilon_{\text{mach}}$ denotes the machine precision, namely,
\[
F^{(k-1)} - F^{(k)} \le \epsilon_{\textsc{bcd}}\, F^{(k-1)},
\]
where $F^{(k)}$ denotes the tour length after cycle~$k$.

The feasible set $B = B(c_1,r_1) \times B(c_2,r_2) \times \dots \times B(c_n,r_n) \subset \mathbb{R}^{2n}$ of~\eqref{probfixedseq} is compact and convex, since it is the Cartesian product of $n$ closed and bounded disks. The total tour length objective function in~\eqref{probfixedseq} is continuous and convex. Under the additional assumption that consecutive disks are disjoint, it is also continuously differentiable on $B$. Assuming that a global minimizer of each single-disk subproblem is computed at every step, the cyclic update scheme is well defined and generates an infinite sequence of iterates with monotonically non-increasing objective values. Since $B$ is compact and the objective function is continuous and convex, it is pseudoconvex on $B$, and every level set of the objective function relative to $B$ is compact. Therefore, the convergence result of the block nonlinear Gauss--Seidel method for pseudoconvex objectives~\citep{grippo2000convergence} implies that every limit point of the sequence is a global minimizer of~\eqref{probfixedseq}. In the numerical experiments section, we will analyze what happens when the disks intersect.

\subsection{Single-disk subproblem}

In this section, we introduce a method for solving the single-disk subproblem~\eqref{singlediskprob}. To simplify the presentation, we assume throughout this section that the fixed visiting sequence is given by $\pi(i)=i$ for $i=1,\dots,n$, with indices interpreted cyclically modulo $n$, so that $n+1$ denotes $1$ and $0$ denotes $n$. At each iteration of the BCD algorithm, the procedure optimizes the location of a single visit point~$y_i$ while holding all other spatial coordinates fixed. Given the fixed sequence and fixed locations of the predecessor $y_{i-1}$ and successor $y_{i+1}$, the single-disk subproblem for customer $i$ seeks a point $y_i \in \mathbb{R}^2$ that minimizes the sum of Euclidean distances to these two fixed adjacent nodes while remaining inside customer $i$'s circular neighborhood. Formally, this subproblem is formulated as the constrained convex optimization problem:
\begin{equation} \label{single-disk}
\begin{array}{rl}
\displaystyle\Minimize_{y_i \in \mathbb{R}^2} & \| y_i - y_{i-1} \|_2 + \| y_i - y_{i+1} \|_2 \\[4mm]
\mbox{subject to} & \| y_i - c_i \|_2 \le r_i.
\end{array}
\end{equation}
Geometrically, this subproblem corresponds to finding the point $y_i$ within the disk $B(c_i,r_i)$ that minimizes the perimeter of the triangle formed by $y_{i-1}$, $y_i$, and $y_{i+1}$. The objective function is convex and continuous, and the feasible region is compact, so a global minimizer exists.

The single-disk continuous subproblem~\eqref{single-disk} can be solved efficiently by exploiting its underlying geometric properties. Because the objective function measures the path length through $y_i$ between two fixed nodes $y_{i-1}$ and $y_{i+1}$, the optimal visit point $y_i^\star$ is guaranteed to lie either on the line segment connecting $y_{i-1}$ and $y_{i+1}$ or on the boundary of $B(c_i, r_i)$. Specifically, if $I:=[y_{i-1},y_{i+1}] \cap B(c_i, r_i) \neq \emptyset$, then any point in $I$ achieves the unconstrained minimum distance $\| y_{i-1} - y_{i+1} \|_2$. Conversely, if $I=\emptyset$, then the optimal point $y_i^\star$ lies on the boundary of $B(c_i, r_i)$, at a point where the angles formed by the two incident segments with the tangent to the circle are equal. In this active-boundary case, the problem reduces to a 1D optimization task over the boundary angle $\theta \in [0, 2\pi)$.

In the case $I \neq \emptyset$, the two obvious options for the choice of $y_i^\star \in I$ would be the orthogonal projection of $c_i$ onto $I$, or the midpoint of $I$. We adopt a third: 
\begin{equation} \label{athird}
y_i^\star = \operatorname*{arg\,max}_{y \in I} \min\bigl(\|y - y_{i-1}\|_2,\ \|y - y_{i+1}\|_2\bigr),
\end{equation}
the point of $I$ that maximizes its distance to the closer of the two neighboring visit points. This choice tends to keep the iterates away from configurations at which the objective function of either~\eqref{probfixedseq} or~\eqref{single-disk} may be non-differentiable, depending on whether consecutive disks are disjoint. To see why this matters, suppose that two consecutive visit points coincide, i.e., $y_i=y_{i+1}$. When BCD optimizes with respect to $y_i$, its current position is optimal, since it lies on the segment $[y_{i-1},y_{i+1}]$ and therefore attains the lower bound $\|y_{i+1}-y_{i-1}\|_2$ given by the triangle inequality. Analogously, when BCD optimizes with respect to $y_{i+1}$, its current position, which coincides with $y_i$, is also optimal. Hence, if each block update kept the current point whenever it is optimal, these two points would never move again: two coinciding consecutive visit points would constitute a fixed point of BCD, which, in general, has no optimality property with respect to problem~\eqref{probfixedseq}. The maximin choice above prevents this situation whenever $I$ contains more than one point, since it moves $y_i$ to the point of $I$ farthest from both neighbors, thereby separating it from $y_{i+1}$.

When the optimal visit point lies on the boundary of the circular neighborhood, translating coordinates so that $c_i$ coincides with the origin simplifies both the analysis and the numerical solution of the subproblem. Under this shift, the neighbor coordinates become $\hat{y}_{i-1} = y_{i-1} - c_i$ and $\hat{y}_{i+1} = y_{i+1} - c_i$, allowing any candidate visit point on the boundary to be parameterized in polar form as $\hat{y}_i(\theta) = r_i (\cos\theta, \sin\theta)^T$, where $\theta \in [0, 2\pi)$. Thus, subproblem~\eqref{single-disk} reduces to the one-dimensional problem
\begin{equation} \label{single-disk-1D}
\Minimize_{\theta \in [0,2\pi)} f(\theta) := \|\hat{y}_i(\theta) - \hat{y}_{i-1}\|_2 + \|\hat{y}_i(\theta) - \hat{y}_{i+1}\|_2.
\end{equation}
Note that, while problem~\eqref{single-disk} is convex, problem~\eqref{single-disk-1D} is not. (Geometrically, the level sets of the two-dimensional objective function in~\eqref{single-disk} are ellipses with foci at $y_{i-1}$ and $y_{i+1}$. When the feasible region contains neither focus and does not intersect the segment joining them---which is precisely the case considered when solving~\eqref{single-disk-1D}---the confocal ellipses are tangent to the circle $\partial B(c_i,r_i)$ at exactly two points, where the objective function attains its minimum and maximum. Consequently, the objective function in~\eqref{single-disk-1D} has a unique minimum and a unique maximum on $[0,2\pi)$. We emphasize, however, that we do not have a formal proof of this property and mention it only as an interesting fact.) Given a solution $\theta^\star$ to~\eqref{single-disk-1D}, a solution to~\eqref{single-disk} is recovered via $y_i^\star = \hat{y}_i(\theta^\star) + c_i$.

\subsection{Multi-start 1D Newton's method}

\begin{figure}[ht!]
\FIGURE
{\resizebox{\textwidth}{!}{
\begin{tabular}{ccc}
\includegraphics[scale=1.0]{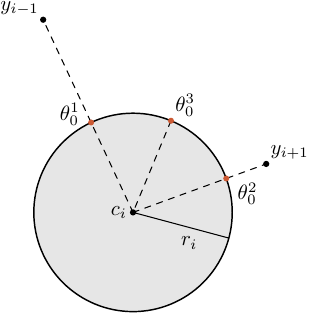} & &
\includegraphics[scale=1.0]{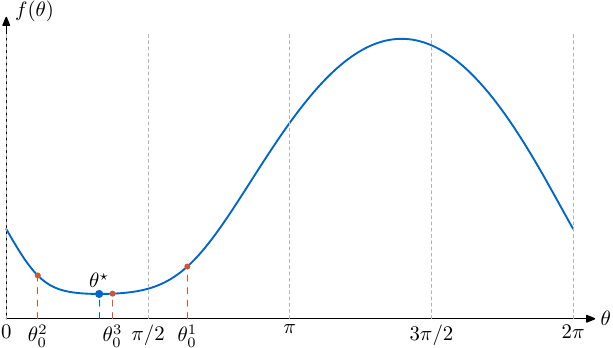}
\end{tabular}}}
{Geometry of the single-disk subproblem and the corresponding univariate function $f(\theta)$. \label{fig1}}
{On the left, geometric arrangement of the disk with center $c_i$ and radius $r_i$ and the points $y_{i-1}$ and $y_{i+1}$. The figure also shows the three anchor angles used as initial guesses in the multi-start minimization process. On the right, we depict the univariate function $f(\theta)$ corresponding to the figure on the left, along with the three anchor angles, and the global minimizer $\theta^\star$. The graph represents the situation that, as we will see later, is the most common in practice: the angles $\theta_0^1$, $\theta_0^2$, and $\theta_0^3$ are close to $\theta^\star$, which causes Newton's method to find the global minimizer when starting from all three of them.}
\end{figure}

As a first approach for solving problem~\eqref{single-disk-1D}, aiming to cover all potential basins of attraction, we consider a multi-start 1D Newton scheme initialized at three fundamental geometric anchor angles:
\begin{enumerate}
\item $\theta_0^1$, the direction of $\hat{y}_{i-1}$: the direct line-of-sight angle from $c_i$ to the predecessor $y_{i-1}$;
\item $\theta_0^2$, the direction of $\hat{y}_{i+1}$: the direct line-of-sight angle from $c_i$ to the successor $y_{i+1}$;
\item $\theta_0^3$, the bisector angle: writing $\hat u_1 = \frac{\hat{y}_{i-1}}{\|\hat{y}_{i-1}\|}$ and $\hat u_2 = \frac{\hat{y}_{i+1}}{\|\hat{y}_{i+1}\|}$, $\theta_0^3$ is the direction of $\hat u_1 + \hat u_2$,
\end{enumerate}
where by the direction of a nonzero vector we mean the angle it makes with the positive $x$-axis. In practice, the direction $\theta$ of a vector $v=(v_1,v_2)$ is given by $\theta=\operatorname{atan2}(v_2,v_1)$ defined as
\[
\operatorname{atan2}(v_2,v_1) = \begin{cases}
\arctan(v_2/v_1) & \text{if } v_1>0,\\[2pt]
\arctan(v_2/v_1)+\pi & \text{if } v_1<0\le v_2,\\[2pt]
\arctan(v_2/v_1)-\pi & \text{if } v_1,v_2<0,\\[2pt]
\pi/2 & \text{if } v_1=0<v_2,\\[2pt]
-\pi/2 & \text{if } v_2<0=v_1.
\end{cases}
\]
We set $\theta^\star = \arg\min_{\ell=1,2,3} f(\theta^{\star,\ell})$, where $\theta^{\star,\ell}$ is the solution found when starting from $\theta_0^\ell$. See Figure~\ref{fig1}. The Newton iteration is given by
\begin{equation}
\theta_{k+1} = \theta_k - \frac{f'(\theta_k)}{f''(\theta_k)},
\label{eq:newton}
\end{equation}
where
\begin{equation} \label{firstder}
f'(\theta) = \frac{ \hat{y}_i'(\theta)^T (\hat{y}_i(\theta) - \hat{y}_{i-1})}{\| \hat{y}_i(\theta) - \hat{y}_{i-1} \|_2} + \frac{\hat{y}_i'(\theta)^T (\hat{y}_i(\theta) - \hat{y}_{i+1})}{\| \hat{y}_i(\theta) - \hat{y}_{i+1} \|_2},
\end{equation}
\[
f''(\theta) = \sum_{j \in \{i-1, i+1\}} \Biggl[ \frac{r_i^2 + (\hat{y}_i(\theta) - \hat{y}_j)^T \hat{y}_i''(\theta)}{\| \hat{y}_i(\theta) - \hat{y}_j \|_2} - \frac{\left( (\hat{y}_i(\theta) - \hat{y}_j)^T \hat{y}_i'(\theta) \right)^2}{\| \hat{y}_i(\theta) - \hat{y}_j \|_2^3} \Biggr],
\]
$\hat{y}_i'(\theta) = r_i (-\sin\theta, \cos\theta)^T$, and $\hat{y}_i''(\theta) = -r_i (\cos\theta, \sin\theta)^T$.

\subsection{Algebraic characterization of the stationary points} \label{algesum}

Besides the numerical one-dimensional optimization approach, the stationary points of problem~\eqref{single-disk-1D} admit an exact algebraic characterization, derived in Appendix~\ref{algeform}. Applying the tangent half-angle substitution $t=\tan(\theta/2)$ to the stationarity condition $f'(\theta)=0$ leads to a quartic polynomial $Q(t)$ with real coefficients. Alternatively, writing $p=\hat y_{i-1}$, $q=\hat y_{i+1}$, and $z=\hat y_i(\theta)=r_ie^{i\theta}$ as complex numbers and invoking the reflection principle of geometric optics, the same condition leads to the quartic polynomial $P(z)=\alpha z^4+\beta z^3+\delta z+\epsilon$, with complex coefficients and no $z^2$ term, all of whose roots lie on the circle $|z|=r_i$. The two polynomials have exactly the same roots under the correspondence $t\leftrightarrow\theta\leftrightarrow z$. Exactly two of the four roots are stationary points of $f$; the other two are stationary points of $h(\theta):=\|\hat y_i(\theta)-\hat y_{i-1}\|-\|\hat y_i(\theta)-\hat y_{i+1}\|$. We have found, empirically, that the two cases are distinguished by the sign of $\operatorname{Re}[(q-z)/(p'(z)-z)]$, where $p'(z)$ is the mirror image of $p$ across the tangent to the circle at $z$: it is negative exactly at the stationary points of $f$, as we verify in Section~\ref{evalthree}. Since the closed-form solution of a quartic by Ferrari's method is impractical and numerically delicate, we compute the roots of $P(z)$ numerically, either by the Durand--Kerner method or as the eigenvalues of its $4\times 4$ companion matrix, computed with LAPACK's zgeev. These two procedures provide an independent means of computing the stationary points of~\eqref{single-disk-1D}, against which the Newton approach is assessed in Section~\ref{sec41}.

\section{Numerical experiments with the BCD and its subproblem} \label{sec4}

In this section, we present two numerical experiments with different purposes. In Section~\ref{sec41}, we test the three proposed approaches for solving the single-disk subproblem~\eqref{singlediskprob} or \eqref{single-disk}. In Section~\ref{sec42}, we assess the performance of the proposed BCD method, equipped with the selected single-disk solver, for computing the optimal visit points associated with a given fixed visiting sequence, i.e.\ BCD applied to problem~\eqref{probfixedseq}. We compare its performance with that of CPLEX applied to the SOCP formulation~\eqref{socp-fixed-sequence} of the same problem. Experiments with BCD within two metaheuristics for the CETSP are presented in Section~\ref{secmeta}.

The one-dimensional Newton method, the two alternatives for computing the roots of the quartic polynomial, the proposed BCD method, and the ILS metaheuristic were implemented in Fortran~90. To facilitate reproducibility, the complete source code is publicly available at \url{https://github.com/rsmatos11/berma2026-cetsp}. To promote broader dissemination and facilitate the potential use of the proposed BCD method combined with the one-dimensional Newton method for computing the visit points of a fixed visiting sequence, a C++ implementation of this procedure is also available. The benchmark instances, detailed solutions, and tour visualizations are provided alongside the source code.

The experiments were carried out on an Intel Core i9-12900K (12th Gen) processor with 16 cores (24 threads) operating at up to 5.2 GHz and 128 GB of RAM, running Ubuntu 24.04.3. The code was compiled using the GNU Fortran compiler GCC (version 13.3.0) with the \texttt{-O3} optimization directive enabled. The SOCP models were solved using IBM ILOG CPLEX Optimization Studio 22.1.2 with all default parameters.

\subsection{Solving the non-trivial case of the single-disk subproblem} \label{sec41}

\subsubsection{Evaluation of the three alternatives} \label{evalthree}

We designed a Monte Carlo experiment to evaluate the numerical methods for solving problem~\eqref{single-disk} in the case in which $[y_{i-1},y_{i+1}] \cap B(c_i,r_i) = \emptyset$. Instances were generated by rejection sampling as follows. The disk is fixed, without loss of generality, at the origin, with radius $r_i$ drawn uniformly in $[R_{\min},R_{\max}]=[0.5,10]$. For each instance, an outer spatial scale $\mathrm{BOX}$ is drawn independently and uniformly in $[\mathrm{BOX}_{\min},\mathrm{BOX}_{\max}]=[1.0,10.0]$. Points $y_{i-1}$ and $y_{i+1}$ are each drawn independently in the annular shell between radius $r_i$ and $\sqrt{8}\,\mathrm{BOX}$ around the origin, with the radial coordinate drawn as $r=\sqrt{r_i^2+u\,(8\,\mathrm{BOX}^2-r_i^2)}$ for $u$ uniform in $[0,1)$ and the angle uniform in $[0,2\pi)$. Instances for which the segment $[y_{i-1},y_{i+1}]$ intersects the disk are discarded; sampling continues until exactly $N=10^6$ instances have been accepted.

The stopping criteria were tailored to each solver's respective characteristics to ensure optimal convergence and an accuracy of the order of machine precision. For the 1D Newton solver, the iterations for each anchor angle were terminated when the absolute difference between successive angular iterates fell below a tolerance, namely, $|\theta^{(k+1)} - \theta^{(k)}| < 10^{-13}$, or if the maximum allowed iterations ($k_{\max} = 50$) was reached. For the Durand--Kerner method, which simultaneously refines all four complex roots of the quartic polynomial, the iterations were stopped when the maximum infinity norm of the updates across all roots satisfied $\max_i |z_i^{(k+1)} - z_i^{(k)}| < 10^{-13}$, also bounded by $k_{\max} = 200$. In contrast, the companion-matrix approach relies on LAPACK's zgeev routine, where the stopping criterion is internal and does not take the form of a main iteration loop which the user could stop by adding a stopping condition: it employs the standard QR algorithm, which terminates when the off-diagonal elements of the upper Hessenberg matrix are reduced to zero within floating-point machine precision ($\epsilon_{\text{mach}} \approx 2.22 \times 10^{-16}$).

In terms of numerical efficacy and computational efficiency, the 1D Newton method demonstrated a clear practical advantage over the polynomial root-finding approaches. Across all test instances, 1D Newton produced the correct global solution in $100.00\%$ of instances, matching the Durand--Kerner and companion-matrix eigensolvers to within $1.4\times10^{-14}$ in every case ($=64\,\epsilon_{\text{mach}}$ for any pair, a floating-point rounding floor rather than a genuine algorithmic discrepancy). This held despite at least one of the three anchor angles individually failing to converge within $k_{\max}=50$ iterations in $962$ instances ($0.0962\%$ of the dataset): in every such case, one of the remaining anchor angles still converged to the correct global minimizer, so the final reported solution was unaffected in all instances. The Durand--Kerner and companion-matrix eigensolvers, for comparison, found exactly four roots of $P(z)$ on the circle $|z|=r_i$ in every single instance, and we used this dataset to verify the sign criterion of Section~\ref{algesum} for distinguishing genuine from spurious roots: in every instance, the roots for which $\operatorname{Re}[(q-z)/(p'(z)-z)]$ is negative were exactly two, and these always coincided exactly with the two stationary points of $f$ identified by checking whether $|f'(\theta)|$ falls below a small tolerance.

Computationally, 1D Newton proved the most efficient, taking an average of $8.93\times10^{-7}$ seconds per call across the full dataset --- averaging $15.3$ Newton iterations per instance, summed over its three anchors --- making it $1.75\times$ faster than the Durand--Kerner polynomial solver ($1.56\times10^{-6}$ seconds per call, averaging $21.8$ Durand--Kerner iterations) and $6.64\times$ faster than the companion-matrix eigenvalue approach via LAPACK's zgeev ($5.93\times10^{-6}$ seconds per call).

\subsubsection{Reducing the number of anchors: an adaptive Newton scheme} \label{secadapt}

Running Newton from all three anchors is usually unnecessary. An analysis of how the anchors' failures overlap in the dataset above, detailed in Appendix~\ref{appadapt}, shows that $\theta_0^1$ and $\theta_0^2$ fail jointly only when $f$ is nearly flat, which is precisely the regime in which the bisector anchor $\theta_0^3$ is reliable, and that starting from whichever of $\theta_0^1$ and $\theta_0^2$ has the smaller value of $f$ almost always suffices. This leads to the following adaptive multi-start procedure: (a) run Newton starting from whichever of $\theta_0^1$ and $\theta_0^2$ has the smaller value of $f$; (b) if $|f(\theta_0^1)-f(\theta_0^2)|/\min(f(\theta_0^1),f(\theta_0^2)) < 10^{-4}$, also run Newton starting from the other one; (c) if $\rho_i < 0.05$, where $\rho_i = r_i\|\hat u_1+\hat u_2\| / (\|\hat y_{i-1}\|+\|\hat y_{i+1}\|)$ measures the flatness of $f$, also run Newton starting from $\theta_0^3$.

Across all the test instances, this adaptive multi-start procedure used, on average, $1.18$ anchors and $6.32$ Newton iterations per instance --- compared to $3$ anchors and $15.28$ iterations for the unconditional scheme that always tries $\theta_0^1$, $\theta_0^2$, and $\theta_0^3$ --- while reproducing the unconditional scheme's global minimizer exactly, with zero disagreements across all instances. This translates into an average solve time of $4.80\times10^{-7}$ seconds per instance, $1.86\times$ faster than the unconditional three-anchor Newton scheme ($8.93\times10^{-7}$ seconds), $3.26\times$ faster than the Durand--Kerner quartic solver ($1.56\times10^{-6}$ seconds), and $12.36\times$ faster than the companion-matrix eigensolver ($5.93\times10^{-6}$ seconds).

\subsection{BCD with multi-start 1D Newton versus SOCP formulation by CPLEX} \label{sec42}

We now compare, on fixed visiting sequences rather than single disks, the proposed BCD-based procedure against CPLEX solving the SOCP formulation~\eqref{socp-fixed-sequence}, for which CPLEX employs a logarithmic-barrier interior-point method. We generated 2,000 fixed-sequence CETSP instances: 200 values of~$n$, ranging from~10 to~408 in increments of~2, with~10 independently seeded random instances per value of~$n$. Disk centers are uniform in $[0,10]^2$ and radii are uniform in $[r_{\min}(n),r_{\max}(n)]$, with $r_{\max}(n)=3/\sqrt n$ and $r_{\min}(n)=0.1\,r_{\max}(n)$. The expected distance from a disk center to the center of its nearest neighbor is about $5/\sqrt n$ for $n$ uniform points in a box of side $10$. These choices keep the ratio between the radii and this spacing between centers, and hence the degree of overlap between the disks, roughly constant as $n$ grows: $r_{\max}(n)$ is about $0.6$ times this spacing, which gives a mix of overlapping and non-overlapping disks. The fixed visiting sequence is simply the disk generation order, i.e.\ an arbitrary sequence. CPLEX was restricted to a single thread, to match BCD's inherently sequential execution.

Using the center of each disk as the initial visit point, BCD converged in a median of 4 sweeps (mean 4.37, range 3--39), while CPLEX required a median of 22 barrier iterations (mean 21.47, range 13--26). In terms of CPU time, BCD is markedly faster, with a mean solve time of $8.95\times10^{-4}$ seconds compared with $1.51\times10^{-2}$ seconds for CPLEX, yielding a mean speedup of $16.8\times$. The comparison of solution quality requires closer examination. Across all 2{,}000 instances, BCD appears to produce longer tours than those reported by CPLEX for the SOCP formulation, with a mean relative difference of $5.7 \times 10^{-4}\%$ (maximum $2.0\times10^{-3}\%$). This apparent gap is almost entirely attributable to residual infeasibility in the CPLEX solutions. Under CPLEX's default tolerances, the mean of the maximum disk constraint violations per instance, measured as the distance from the visit point to its disk (which should be zero in exact math, as the visit point must belong to the disk), is $1.13\times10^{-4}$. BCD, by construction, satisfies every disk constraint to machine precision, with a mean violation of $8.78\times10^{-16}$, approximately eleven orders of magnitude smaller than those observed for CPLEX. Tightening CPLEX's QCP convergence and feasibility tolerances reduces both the mean objective gap relative to BCD and the mean disk constraint violation, but increases solve time and leaves some instances unsolved. Alternatively, projecting the visit points computed by CPLEX onto their disks yields feasible solutions and essentially eliminates the gap: the mean relative difference between the tour lengths of BCD and of the projected CPLEX solutions becomes $-7\times10^{-7}\%$. The negative sign indicates that, on average, the BCD tours are slightly shorter than the projected CPLEX tours; this occurs in 1,861 of the 2,000 instances, although the differences are negligible. The relative difference exceeds $10^{-5}\%$ in only 5 of the 2,000 instances, with a maximum of $3.7\times10^{-4}\%$. These results show that BCD and CPLEX compute solutions of essentially equivalent optimality accuracy, while BCD achieves substantially greater feasibility accuracy and is much faster on the tested instances.

The choice~\eqref{athird} of the point of $I$ when $I\neq\emptyset$ has a marked effect on these results. With the maximin rule adopted here, the tour length obtained by BCD exceeds that of the projected CPLEX solution by more than $10^{-5}\%$ in only 5 instances, with a maximum of $3.7\times10^{-4}\%$. Since the projected CPLEX solutions are feasible, these few cases correspond to BCD stopping at points that are not global minimizers, albeit with negligible differences. If the midpoint of $I$ is used instead, this happens in 25 instances, with a maximum of $2.3\times10^{-3}\%$, and if the orthogonal projection of $c_i$ onto $I$ is used, it happens in 325 instances, with differences above $0.1\%$ in 22 of them and a maximum of $2.0\%$. This confirms that BCD can stall at non-optimal fixed points when consecutive visit points coincide, and that the maximin rule is an effective safeguard against this behavior. 

Parallelization does not improve the performance of either solver at this instance scale. Allowing CPLEX to use all available cores instead of one increases its mean solve time from $1.51\times10^{-2}$ to $2.44\times10^{-2}$ seconds. Similarly, an OpenMP variant of BCD that updates the even-indexed visit points concurrently, followed by the odd-indexed ones, is slower than the sequential version ($2.66\times10^{-3}$ versus $8.95\times10^{-4}$ seconds). Even so, parallel BCD remains approximately $9.2\times$ faster than parallel CPLEX. This pattern persists on newly generated instances with up to $n=50{,}000$, more than 120 times the size of the largest instance considered above ($n=408$). Across these tests, the ratio of parallel to sequential solve times ranges from $1.13$ to $1.66$ for BCD and from $1.13$ to $2.54$ for CPLEX, with neither solver benefiting from parallelization and CPLEX's slowdown becoming more pronounced at the largest sizes tested. 

A few remarks on the number of BCD sweeps are in order. The largest sweep counts occur on instances where two adjacent visit points end up close to each other without coinciding. Among the 282 instances (14.1\%) in which the closest adjacent pair is separated by less than 0.05, BCD requires a mean of 8.07 sweeps (maximum 39), compared with 3.77 sweeps (maximum 16) for the remaining 1,718 instances (85.9\%). Thus, the mean sweep count more than doubles in cases where the theory does not guarantee global convergence. Nevertheless, the numerical results consistently indicate convergence to a global minimizer, albeit with more sweeps.

\section[Numerical experiments of BCD within trajectory-based and population-based metaheuristics]{Numerical experiments of BCD within trajectory-based and\\ population-based metaheuristics} \label{secmeta}

In this section, we assess BCD within two metaheuristics for the CETSP. In Section~\ref{sec43}, we consider a trajectory-based metaheuristic, an Iterated Local Search (ILS), in which BCD computes the visit points of the candidate visiting sequences. In Section~\ref{secma}, we consider a population-based metaheuristic, the memetic algorithm MA-CETSP of~\citet{Lei2024}, one of the best-performing methods for the CETSP, in which the general-purpose SOCP solver is replaced by BCD, among other modifications.\par

\subsection{ILS-BCD: Iterated Local Search with BCD} \label{sec43}

In this section, to illustrate the use of BCD for computing visit points within a trajectory-based metaheuristic, we consider the Iterated Local Search (ILS) metaheuristic applied in~\citet{ber2025-tspn} to the TSPN with polygonal neighborhoods. The method combines a search over visiting sequences with the computation of visit points by BCD. Its adaptation to the CETSP retains the combinatorial search structure and uses the proposed BCD method with multi-start 1D Newton to compute visit points within the disks.

A solution is represented by a permutation $\pi$ of $\{1,\dots,n\}$, that stands for the visiting sequence, plus the visit points $y_{\pi(1)},\dots,y_{\pi(n)}$, satisfying $y_{\pi(i)} \in B(c_{\pi(i)},r_{\pi(i)})$ for $i=1,\dots,n$. We denote its tour length by
$$
F(\pi,y) = \sum_{j=1}^{n}\|y_{\pi(j)}-y_{\pi(j-1)}\|_2.
$$
For a fixed sequence $\pi$, minimizing $F(\pi,y)$ over $y_i \in B(c_i,r_i)$ for $i=1,\dots,n$ is precisely the continuous subproblem considered in Section~\ref{sec3}.

The initial solution fed to ILS is constructed using the nearest-neighbor heuristic applied to the disk centers $c_1,\dots,c_n$. The heuristic is run from each customer, and BCD is applied to each resulting sequence, starting from $y_i=c_i$ for $i=1,\dots,n$. The shortest resulting tour is selected and refined by local search. Each subsequent ILS iteration consists of perturbing the current solution, applying local search, and deciding whether to accept the resulting solution as the new current solution. The best solution encountered throughout the search is retained. Local search uses Variable Neighborhood Descent (VND). The selected neighborhoods are explored sequentially, and the search restarts from the first whenever an improving solution is accepted. It terminates when none of the neighborhoods yields an improvement. Each move modifies the visiting sequence $\pi$, after which BCD, limited to at most $s_{\max}$ sweeps, recomputes the visit points $y$ to evaluate the resulting tour. Optionally, once VND terminates, a complete BCD (up to 200 sweeps) is applied to the solution it returns. The VND uses eight neighborhood structures (2-opt, 3-opt, double bridge, swap, insertion, string relocation, Balas--Simonetti, and Gutin), described in Appendix~\ref{appvnd}.

The perturbation applies $\ell$ consecutive double-bridge moves, with the edges involved in each move selected at random. Each move acts on the sequence produced by the preceding one. After the complete perturbation, BCD, limited to at most $s_{\max}$ sweeps, recomputes the visit points, and VND refines the resulting solution. The parameter $\ell$ therefore controls the perturbation strength.

A solution that improves upon the current solution is always accepted. A solution of higher cost may also be accepted if
$$
F(\bar \pi,\bar y) \leq F(\pi^\star,y^\star) + \eta\,|F(\pi^\star,y^\star)|,
$$
where $(\bar \pi,\bar y)$ is the solution returned by local search and $(\pi^\star,y^\star)$ is the best solution found so far. The threshold $\eta$ starts at $\eta_0>0$ and is progressively reduced according to $\eta\leftarrow h\eta$, where $h\in(0,1)$. Following the calibration convention in~\citet{ber2025-tspn}, $\kappa\geq 1$ denotes the number of threshold reductions over the run. The iterations at which the reductions occur are determined from the time limit: the first iterations of the ILS are used to estimate how long an iteration takes, this estimate gives how many iterations still fit in the remaining time, and the reductions are spread evenly over them. The parameters of ILS-BCD were calibrated with \texttt{Irace}~\citep{LopezIbanez2016}; the calibration procedure and the selected values are described in Appendix~\ref{seccal}.

\subsubsection{ILS-BCD versus ILS-CPLEX}

We compared two versions of ILS that differ only in the method used to compute visit points for a given visiting sequence: the proposed BCD method, and CPLEX for solving the SOCP formulation~\eqref{socp-fixed-sequence}. We refer to them as ILS-BCD and ILS-CPLEX, respectively. For ILS, we adopted the parameters specified in Appendix~\ref{seccal}, except that $s_{\max}$ was set to~$+\infty$ instead of~1. This ensures that the subproblem of determining the visit points for a given sequence is effectively solved by BCD, just as it is by CPLEX.  The comparison used 120 instances from Section~\ref{sec42}: 60 with $n\in\{10,12,14,16,18,20\}$ and 60 with $n\in\{30,34,38,42,46,50\}$, with ten instances per size. Each instance was solved once by each version, using the same random seed and a time limit of 3,600~seconds. Table~\ref{tab:ils-bcd-vs-cplex} summarizes some details. For each size group, it reports the mean relative difference between the tour lengths obtained by the two versions and, for each version, the mean CPU time, number of ILS iterations, and number of calls to the visit-point solver up to the moment at which the best solution was found, as well as the CPU time per call. Both versions produced essentially identical tour lengths, with mean relative differences of 0.00014\% and 0.00027\% across the two size groups, respectively,  attributable to CPLEX's feasibility tolerances.

Specifically, CPLEX solves the SOCPs using a barrier method, a primal-dual interior-point algorithm, which returns an approximate solution that may slightly violate the disk constraints (by $2.3\times10^{-5}$ on average for these 120 instances in the experiments of Section~\ref{sec42}). This minor violation allows the visit points to lie slightly outside the disks, thereby yielding slightly shorter tours. In contrast, BCD yields solutions where the visit points lie within the disks up to machine precision, i.e., on the order of $10^{-16}$. ILS-BCD required substantially less CPU time to reach its final recorded incumbent: on average, 0.016 versus 0.884~seconds in the smaller group and 4.43 versus 390.84~seconds in the larger group. The corresponding ratios of mean times are $55.25$ and $88.23$, exceeding the speedup of $16.8$ observed in Section~\ref{sec42} for fixed-sequence evaluations initialized at the disk centers.

\begin{table}[ht!]
\TABLE
{Comparison between ILS-BCD and ILS-CPLEX on 120 instances, grouped by size. \label{tab:ils-bcd-vs-cplex}}
{  \small
  \begin{tabular}{crrrrrrrrrr}
  \hline
  & & \multicolumn{2}{c}{CPU time (s)}
  & \multicolumn{2}{c}{ILS iterations}
  & \multicolumn{2}{c}{Calls}
  & \multicolumn{2}{c}{Time per call (ms)}
  & \\
  \cmidrule(lr){3-4} \cmidrule(lr){5-6} \cmidrule(lr){7-8} \cmidrule(lr){9-10}
  $n$ & Gap (\%) & BCD & CPLEX & BCD & CPLEX & BCD & CPLEX & BCD & CPLEX & Ratio \\
  \hline
  10--20 & 0.00014 & 0.016 & 0.884   & 0.1 & 0.6  & 227    & 692     & 0.070 & 1.28 & 18 \\
  30--50 & 0.00027 & 4.43  & 390.84  & 9.8 & 23.5 & 44,949 & 116,518 & 0.099 & 3.35 & 34 \\
  \hline
  \end{tabular}}
{Gap is the mean relative difference between the tour lengths obtained with BCD and with CPLEX. CPU time, ILS iterations, and calls are means measured up to the moment at which the best solution was found; time per call is the CPU time divided by the number of calls to the visit-point solver; Ratio is the ratio between the times per call of CPLEX and BCD.}
\end{table}

Table~\ref{tab:ils-bcd-vs-cplex} shows that this greater advantage results from two effects, whose product approximately gives the ratios of mean times reported above ($3.0\times 18\approx 55$ and $2.6\times 34\approx 88$). First, ILS-CPLEX required more ILS iterations and about 3.0 and 2.6 times as many calls to the visit-point solver to reach its final incumbent in the smaller and larger groups, respectively (692 versus 227 and 116,518 versus 44,949 calls). This is explained by the feasibility tolerance of CPLEX. Its residual disk constraint violations can produce small variations in the reported tour length when evaluating identical or equivalent sequences, and these variations may be accepted as improvements, causing additional incumbent updates and delaying the final recorded incumbent without a meaningful reduction in tour length. Since BCD maintains disk feasibility to approximately machine precision, ILS-BCD recognizes more reliably when a candidate sequence offers no further improvement. Second, the CPU time per call was 18 and 34 times lower with BCD in the two groups, respectively. In the smaller group, this ratio is close to the speedup of $16.8$ observed in Section~\ref{sec42} for fixed-sequence evaluations initialized at the disk centers, whereas in the larger group it is about twice as large. We attribute this difference to the warm start exploited by BCD: after a local-search move modifies the visiting sequence, BCD starts from the visit points of the sequence from which the new one was obtained, rather than from the disk centers, and, since neighboring sequences often differ only slightly, many of these points are already close to optimal for the new sequence. The barrier method used by CPLEX does not exploit the warm start provided by the inherited visit points as effectively as our method does. Thus, the larger speedups reflect both the avoidance of numerically induced incumbent updates and the effective reuse of visit points by BCD, making BCD more advantageous within ILS than the isolated fixed-sequence comparison suggests.

\subsubsection{ILS-BCD on small benchmark instances}

In this section, we assess the solution quality and computational effort of ILS-BCD on the 720 small benchmark instances used in~\citet{Coutinho2016}, for which optimal values are known. Originally proposed in~\citet{Behdani2014}, the set comprises 240 configurations of points randomly placed in a rectangle of dimensions $16\times10$. There are 30 configurations for each of eight sizes, with $7,9,11,\dots,21$ customers, including the depot. Each configuration is used with three common radii for the customer disks, $r\in\{0.25,0.50,1.00\}$, yielding 240 instances per radius and 720 instances in total. This set allows us to compare the tour lengths obtained by ILS-BCD with the optimal values and to examine how the radius affects the computational effort.

We applied ILS-BCD to all 720 instances and compared the resulting tour lengths with these reference values. For ILS, we adopted the parameters specified in Appendix~\ref{seccal}, except that $s_{\max}$ was set to~$+\infty$ instead of~1; this setting was chosen because the instances are small. ILS-BCD matched the reported optimal value in every instance. Table~\ref{tab:ils-small-benchmarks} summarizes the results by radius, reporting the number of instances for which the optimum was attained, the mean and maximum CPU times, the mean and maximum numbers of ILS iterations and BCD calls, and the mean number of BCD sweeps per call. Mean CPU times ranged from 0.016 to 0.048~seconds, and no instance required more than 5.44~seconds. In approximately $90\%$ of the instances, an optimal solution was obtained during initialization, before any perturbation-based ILS iteration. These timings characterize the effort required by the heuristic; they are not directly comparable with the time required by an exact branch-and-bound algorithm to certify optimality.

\begin{table}[ht!]
\TABLE
{Results of ILS-BCD on the 720 small benchmark instances, grouped by the common radius of the customer disks. \label{tab:ils-small-benchmarks}}
{\small
\begin{tabular}{crrrrrrrr}
\hline
& & \multicolumn{2}{c}{CPU time (s)}
& \multicolumn{2}{c}{ILS iterations}
& \multicolumn{2}{c}{BCD calls}
& BCD sweeps \\
\cmidrule(lr){3-4} \cmidrule(lr){5-6} \cmidrule(lr){7-8}
Radius & Optima/Total & Mean & Max. & Mean & Max. & Mean & Max. & per call \\
\hline
0.25 & 240/240 & 0.016 & 0.49 & 0.34 & 12 & 413 & 12,226 & 4.77 \\
0.50 & 240/240 & 0.048 & 5.44 & 0.55 & 61 & 702 & 60,462 & 7.03 \\
1.00 & 240/240 & 0.033 & 1.14 & 0.39 & 22 & 589 & 19,695 & 11.44 \\
\hline
\end{tabular}}
{The last column reports the mean number of BCD sweeps per call.}
\end{table}

The mean number of BCD sweeps per call increases from 4.77 to 7.03 and then to 11.44 as the radius grows. This behavior is consistent with the observations in Section~\ref{sec42}: larger disks allow more overlap, increasing the likelihood that adjacent visit points become very close without coinciding, a situation in which BCD may converge more slowly. This increase in the effort per BCD call does not translate into a monotonic increase in the overall CPU time, which also depends on the number of calls made during the search.

\subsubsection{ILS-BCD on large benchmark instances}

In this section, we assess the solution quality and computational effort of ILS-BCD on the 62 large benchmark instances introduced in~\citet{mennell2009thesis}, for which the reference values are the tour lengths obtained by the heuristics evaluated in that work. For an instance in which all disks have the same radius, the overlap ratio is defined as the ratio between the radius of the disks and the side of the smallest square containing all disks~\citep{mennell2009thesis}. Two groups of instances are considered. The first consists of the 48 instances, each with identical disks, reported in~\citet[Table~2.2]{mennell2009thesis}: 21 instances derived from seven TSPLIB instances (kroD100, rat195, lin318, rd400, pcb442, d493, and dsj1000), each combined with three overlap ratios (0.02, 0.10, and 0.30), and 27 instances with varied overlap ratios (team1-100 to team6-500, bonus1000, rotatingDiamonds1 to rotatingDiamonds5, bubbles1 to bubbles9, concentricCircles1 to concentricCircles5, and chaoSingleDep). The second group consists of the 14 instances with variable radii reported in~\citet[Table~C1]{mennell2009thesis}, built on the same seven TSPLIB instances and on the instances team1-100 to team6-500 and bonus1000. All instances include a depot, which, in our instances, corresponds to an extra customer whose disk has zero radius. The number $n$ of customers, including the depot, ranges from 17 to 1001.

To keep the total running time of the experiment within reasonable limits, each instance was solved once by ILS-BCD, using the parameter configuration selected in the calibration described in Appendix~\ref{seccal}, with a CPU time limit of 10~hours. For instances with $n \leq 200$, we additionally terminated the search after 5,000 consecutive non-improving iterations. Twenty of the 22 instances with $n\le 200$ terminated by this stagnation criterion. Tables~\ref{tab_ils_effort} and~\ref{tab_ils_effort_var} in Appendix~\ref{appresults} report the results for the instances with a common radius and with variable radii, respectively. The first columns identify each instance: its name, the number $n$ of customers, and its overlap ratio (column OR), which is not reported for the instances with varied overlap ratios or with variable radii. For each instance, the table also reports the tour length obtained by ILS-BCD and the computational effort measured up to the moment at which the best solution was found: the CPU time (in seconds), the number of ILS iterations, and the number of BCD calls. On average, the best solution was found after 1.1~hours for the instances with $n\le 200$ (median 0.3~hours) and after 7.9~hours for the instances with $n>200$ (median 9.0~hours). Figure~\ref{fig:berma2026avgtimeRUN5} shows the evolution over time of the average incumbent tour length, normalized by the final value of each instance. For the instances with $n\le 200$, the average incumbent was within 0.2\% of the final value after 10~minutes and within 0.1\% after one hour, and practically no improvement occurred after five hours. For the instances with $n>200$, the average incumbent was 4.4\% above the final value after 30~minutes, 1.9\% after two hours, and 0.6\% after five hours, and it was still decreasing at the end of the run: half of these instances found their best solution after 9.0~hours. This indicates that, for the largest instances, longer runs could further improve the results.

\begin{figure}[ht!]
\FIGURE
{\begin{subfigure}[b]{0.49\textwidth}
\input{berma2026fig2a.tex}
\caption{Instances with $n \le 200$.}
\label{fig:avgtime_a}
\end{subfigure}
\hfill
\begin{subfigure}[b]{0.49\textwidth}
\input{berma2026fig2b.tex}
\caption{Instances with $n > 200$.}
\label{fig:avgtime_b}
\end{subfigure}}
{Evolution of the average incumbent solution over time for ILS-BCD with the parameter configuration selected by Irace, over a run of 10 hours. \label{fig:berma2026avgtimeRUN5}}
{To give equal weight to all instances in the computation of the average, the incumbent tour length $\mathrm{TD}(t)$ of each instance at time $t$ is normalized by the best solution found for that instance at the end of the run, $\mathrm{TD}_{\mathrm{best}}$. The average is taken over the 22 instances with $n \le 200$ in (a) and over the 40 instances with $n > 200$ in (b).}
\end{figure}

Tables~\ref{tab_cetsp2.1} and~\ref{tab_cetsp4.1} in Appendix~\ref{appresults} compare the tour lengths obtained by ILS-BCD with those reported in~\citet{mennell2009thesis} for each heuristic. Since the reference values are reported in~\citet{mennell2009thesis} with two decimal places, all values were rounded to two decimal places for the comparison. Column ``Best'' gives, for each instance, the smallest value among ILS-BCD and the heuristics. Row ``Wins'' counts the instances in which each method attains this value, ties included; row ``Avg Gap'' gives the mean relative difference, in percent, with respect to ``Best''; and row ``Avg Gap (Losses)'' gives the same mean restricted to the instances in which the method does not attain ``Best''. Computational times are not compared, since the experiments in~\citet{mennell2009thesis} were carried out in a different computational environment and under different time limits. We note, however, that some of the heuristics evaluated in~\citet{mennell2009thesis} were run for very long times, which, converted to the machine used in this work, would correspond to roughly 13 to 32 days per instance in the most extreme case, whereas each ILS run was limited to 10~hours. On the 48 instances with a common radius, ILS-BCD attained the best value in 44 instances: in 36 of them it was strictly better than all 15 heuristics, and in the remaining 8 it tied with the best of them. Its average gap was 0.13\%, compared with 5.01\% for HYBRID2, the heuristic with the smallest average gap in~\citet{mennell2009thesis}. Relative to the best heuristic value of each instance, ILS-BCD reduced the tour length by 2.18\% on average and by up to 15.05\% (instance bubbles6). The largest average reductions occurred for the instances with overlap ratio 0.10 (3.40\%) and for the 27 instances with varied overlap ratios (2.85\%). ILS-BCD did not attain the best value in four instances: dsj1000 with overlap ratio 0.02 (3.59\% above GTSP2), chaoSingleDep (1.64\% above YUAN+ and LK-SOCP), team5-499 (0.80\% above GTSP2), and dsj1000 with overlap ratio 0.30 (0.01\% above LK-SOCP). All four instances have $n>200$, a range in which ILS-BCD was still improving at the end of the run (Figure~\ref{fig:berma2026avgtimeRUN5}), which suggests that longer runs could further reduce the gap. On the 14 instances with variable radii, ILS-BCD attained the best value in all instances, being strictly better than all 10 heuristics in 11 of them and tying with the best of them in the remaining 3. Relative to the best heuristic value of each instance, it reduced the tour length by 1.76\% on average and by up to 6.53\% (instance pcb442).

Tables~\ref{tab_recent_common} and~\ref{tab_recent_var} in Appendix~\ref{appresults} extend the comparison to four more recent heuristics evaluated on the same 62 instances: SZVNS~\citep{Wang2019}, (lb/ub)Alg~\citep{Carrabs2020}, GA~\citep{diplacido2022cetsp}, and MA-CETSP~\citep{Lei2024}. For MA-CETSP, both the best and the average results of 20 runs are reported. Over the 62 instances, ILS-BCD obtained better solutions than SZVNS and (lb/ub)Alg in 45 and 35 instances, and worse ones in 5 and 13, respectively. Its results are comparable to those of GA, with 25 better and 13 worse solutions. On the 23 instances whose optimal values are known~\citep{Coutinho2016}, ILS-BCD attained the optimum in 19 and remained within 0.08\% of it in the other four. MA-CETSP, one of the best-performing methods for these instances, remains superior: ILS-BCD improved its best result in three instances, two of which (bubbles1 and bubbles3) differ just in the second decimal place due to rounding, and was worse in 29, with an average gap of 0.78\% over the 62 instances. With respect to the average results of MA-CETSP, ILS-BCD was better in 21 instances and worse in 23. These comparisons should be interpreted in light of computing times, since each ILS run was allowed up to 10~hours, whereas the reference methods report running times of seconds to minutes. We emphasize that ILS-BCD is used here only to illustrate the use of BCD within a metaheuristic, and no attempt was made to make it competitive with the state of the art. The strength of MA-CETSP lies mainly in its combinatorial search, which combines a population-based genetic framework with the LKH heuristic for the TSP~\citep{Helsgaun2000} to determine the visiting sequence. The visit points, in contrast, are computed by a greedy approximation during the search and refined by solving an SOCP with a general-purpose solver at the end of each local search~\citep{Lei2024}. The two contributions are therefore complementary. Section~\ref{secma} combines them: it replaces the general-purpose SOCP solver of MA-CETSP with BCD and evaluates the resulting method on the same 62 instances.

\subsection{MA-CETSP-BCD: Memetic Algorithm with BCD} \label{secma}

The memetic algorithm MA-CETSP of~\citet{Lei2024} is, to the best of our knowledge, one of the best-performing methods for the benchmark instances of~\citet{mennell2009thesis}. As discussed in Section~\ref{sec43}, its strength lies in its combinatorial search, whereas the visit points are computed by a greedy approximation and by an SOCP solved with a general-purpose solver. In this section, we replace this solver with BCD and use the speed and accuracy of BCD and of the single-disk solver of Section~\ref{sec3} to modify other steps of the algorithm (described below). We call the resulting method MA-CETSP-BCD.

\subsubsection{The MA-CETSP}

MA-CETSP maintains a population of 20 solutions, whose initial members are built from k-means clusterings of the disk centers. At each iteration, two parents are chosen at random and combined by a crossover operator into a single offspring. The visit points of the offspring are adjusted by a greedy procedure that moves each visit point according to the visit points of its two neighbors in the tour; the visiting sequence is then improved by the LKH heuristic~\citep{Helsgaun2000}, and the greedy procedure is applied again. Next, a local search with relocation and swap moves (best improvement) is applied, in which the visit point of each moved customer is approximated geometrically, and the visit points of the resulting sequence are finally computed by solving the SOCP~\eqref{socp-fixed-sequence} with Gurobi. The offspring enters the population if it differs sufficiently from the other individuals or improves the best solution, and, when the population exceeds its maximum size, the individuals with the worst fitness, which combines tour length and distance to the rest of the population, are removed. The search stops after 5,000 iterations or after 500 consecutive iterations without improving the best solution.

\subsubsection{Modifications}

MA-CETSP-BCD differs from MA-CETSP in four points.
\begin{enumerate}[leftmargin=*, topsep=2pt, itemsep=2pt, parsep=2pt]
\item \emph{SOCP solver.} Each SOCP is solved by BCD with the adaptive multi-start 1D Newton scheme of Section~\ref{secadapt}, starting from the current visit points and stopping when a sweep reduces the tour length by less than $10^{-6}$ (relative) or after 200 sweeps. At the end of the search, the visit points of the best solution are refined by BCD. Since, in instances with highly overlapping disks, BCD may reduce the tour length very slowly but steadily, the refining stops only when 1,000 consecutive sweeps together reduce the tour length by less than $10^{-6}$ (relative); stopping when a single sweep yields a small reduction left the tour length up to 0.009 above the optimal value for the same visiting sequence in some of these instances.
\item \emph{Exact local search.} In the relocation and swap moves, the visit point of each moved customer is the exact solution of the single-disk subproblem~\eqref{single-disk} for its new neighbors, computed by the same solver, instead of a geometric approximation; a cheap lower bound on the cost of each insertion discards the candidates that cannot improve on the best move found so far.
\item \emph{BCD after LKH.} The greedy procedure applied after LKH is replaced by two BCD sweeps.
\item \emph{Greedy procedure before LKH.} When the segment joining the two neighbors does not intersect the disk, the greedy procedure computes a minimizer of the single-disk subproblem by bisection; we compute it with the 1D Newton method of Section~\ref{sec3} instead. When the segment crosses the disk, we choose the endpoint of the chord closest to one of the neighbors.
\end{enumerate}

\subsubsection{Experimental setup}

Both MA-CETSP and MA-CETSP-BCD were run with the source code of~\citet{Lei2024}, written in C++, in which we implemented the modifications above using the C++ implementation of BCD mentioned in Section~\ref{sec4}. The code was compiled with GCC~13.3.0 (\texttt{-O3 -flto}), and MA-CETSP used Gurobi~12.0.2 restricted to a single thread, the solver of its original implementation (the experiments of Sections~\ref{sec42} and~\ref{sec43} used CPLEX). Each of the 62 instances of~\citet{mennell2009thesis} was solved 20 times by each method, with the same 20 random seeds and the parameters of~\citet{Lei2024}, on the machine described in Section~\ref{sec4}, one run per core; the time limit of 7,200~seconds was never reached. All 1,240 solutions of MA-CETSP-BCD were verified independently of the code: every disk of the instance is visited, every visit point lies in its disk up to $10^{-13}$, and the tour length recomputed from the visit points agrees with the reported one up to $10^{-11}$. The solutions of MA-CETSP violate the disk constraints by at most $1.7\times10^{-5}$, a consequence of the solver's tolerances; projecting their visit points onto the disks changes the reported tour length, rounded to two decimal places, in only 5 of the 1,240 runs. As in~\citet{Lei2024}, tour lengths are compared with two decimal places.

\subsubsection{Results}

Table~\ref{tab_ma_summary} compares MA-CETSP-BCD with MA-CETSP. Considering the best of the 20 runs, MA-CETSP-BCD is better in 20 instances and equal in the remaining 42. Considering the average of the 20 runs, it is better in 37 instances, equal in 23, and worse in 2, with an average relative difference of $-0.26\%$; the Wilcoxon signed-rank test on the 62 averages gives $p=1.4\times10^{-7}$, and the test applied to each instance, with runs paired by seed, indicates a significant improvement (at the 5\% level) in 28 instances and a significant deterioration in 2. The improvements concentrate on the instances with small or varied overlap ratios and, to a lesser extent, on those with variable radii; on the instances with overlap ratio 0.30, both methods obtain the same best values.

\begin{table}[ht!]
\TABLE
{MA-CETSP-BCD versus the original MA-CETSP (SOCP solved by Gurobi), both run with the same 20 seeds on each of the 62 instances of~\citet{mennell2009thesis}. \label{tab_ma_summary}}
{\begin{tabular}{lcccc}
\hline
\multirow{2}{*}{Group} & \multirow{2}{*}{Instances} & Best of 20 & Average of 20 & Avg.\ diff. \\
& & W/T/L & W/T/L & (\%) \\
\hline
\hline
Common radius, overlap ratio 0.02 & 7 & 6/1/0 & 7/0/0 & -0.28 \\
\rowcolor{gray!15}
Common radius, overlap ratio 0.10 & 7 & 2/5/0 & 7/0/0 & -0.53 \\
Common radius, overlap ratio 0.30 & 7 & 0/7/0 & 0/6/1 & +0.00 \\
\rowcolor{gray!15}
Common radius, varied overlap ratios & 27 & 9/18/0 & 17/10/0 & -0.36 \\
Variable radii & 14 & 3/11/0 & 6/7/1 & -0.05 \\
\hline
Total & 62 & 20/42/0 & 37/23/2 & -0.26 \\
\hline
\end{tabular}}
{W/T/L: number of instances in which MA-CETSP-BCD obtains a better (win), equal (tie), or worse (loss) value, considering the best and the average of the 20 tour lengths rounded to two decimal places. Avg.\ diff.: mean, over the instances of the group, of the relative difference between the average tour lengths (negative values favor MA-CETSP-BCD).}
\end{table}

Table~\ref{tab_ma_instances} in Appendix~\ref{appresults} reports the results for each instance, together with the best values reported in~\citet{Lei2024}. Our runs of MA-CETSP reproduce these values closely (equal in 42 instances, better in 8, and worse in 12). MA-CETSP-BCD improves on the best reported value in 21 instances and matches it in the other 41, and it attains the optimal value in the 23 instances whose optimal values are known~\citep{Coutinho2016}. These are 21 instances rather than the 20 of Table~\ref{tab_ma_summary} because the reference here is the set of values reported in~\citet{Lei2024}, not our own runs of MA-CETSP. The best tours of MA-CETSP-BCD are not those of MA-CETSP with improved visit points: in none of the 62 instances do the two methods find the same visiting sequence, and in the 20 instances in which MA-CETSP-BCD is better, the two tours differ geometrically (the largest distance between them, as curves in the plane, has median 0.8\% of the tour length). The more accurate evaluation of candidate sequences thus leads the search to better visiting sequences. On the 720 small instances of Section~\ref{sec43}, MA-CETSP-BCD attains the optimal value in all instances, as ILS-BCD does, and in 97.6\% of the individual runs. (For the smallest of these instances, all the solutions generated by MA-CETSP may be so close to each other that the population collapses to a single individual, and the selection of two distinct parents never terminates. In this case, MA-CETSP-BCD obtains the offspring by random swaps applied to the single individual.)

To identify the source of these improvements, we carried out the ablation study summarized in Table~\ref{tab_ma_ablation}, in which the four modifications are switched on one at a time, in the order listed above. Replacing Gurobi with BCD alone (V1) yields a method equivalent to MA-CETSP: it wins and loses in a similar number of instances ($p=0.33$) and its total running time is 75\% of that of MA-CETSP. Thus, BCD replaces the general-purpose solver without loss. The exact local search (V2) brings the largest improvement, and the BCD sweeps after LKH (V3) complete it: with both, MA-CETSP-BCD matches or improves the best of 20 runs of MA-CETSP in all instances. The Newton-based greedy procedure (V4) does not change the solution quality (average of 20 runs: 12 wins, 36 ties, and 14 losses with respect to V3, $p=0.32$), but reduces the total running time by 16\%.

\begin{table}[ht!]
\TABLE
{Ablation study on the 62 instances of~\citet{mennell2009thesis} (20 runs each, same seeds). \label{tab_ma_ablation}}
{\resizebox{\fullwidth}{!}{
\begin{tabular}{lccccccc}
\hline
\multirow{2}{*}{Variant} & SOCP & Local & BCD sweeps & Greedy step & Best of 20 & Average of 20 & Wilcoxon \\
& solver & search & after LKH & before LKH & W/T/L & W/T/L & $p$ \\
\hline
\hline
V0 (MA-CETSP) & Gurobi & approx. & no & original & -- & -- & -- \\
\rowcolor{gray!15}
V1 & BCD & approx. & no & original & 13/37/12 & 21/20/21 & 0.33 \\
V2 & BCD & exact & no & original & 16/41/5 & 33/23/6 & $2.1\times10^{-5}$ \\
\rowcolor{gray!15}
V3 & BCD & exact & yes & original & 20/42/0 & 38/22/2 & $7.0\times10^{-8}$ \\
V4 (MA-CETSP-BCD) & BCD & exact & yes & Newton & 20/42/0 & 37/23/2 & $1.4\times10^{-7}$ \\
\hline
\end{tabular}}}
{Starting from the original MA-CETSP (V0), each variant switches on one modification: the SOCP solver replaced by BCD with the final refining (V1), the exact local search (V2), two BCD sweeps after LKH (V3), and the Newton-based greedy step before LKH (V4, i.e., MA-CETSP-BCD). W/T/L: wins, ties, and losses with respect to V0, considering the best and the average of the 20 tour lengths rounded to two decimal places; $p$: Wilcoxon signed-rank test on the 62 average tour lengths with respect to V0.}
\end{table}

Table~\ref{tab_ma_effort} reports the computational effort. Because the exact local search keeps improving the best solution for longer, the runs of MA-CETSP-BCD are longer than those of MA-CETSP, and so are the times to the best solution of each run. A fairer comparison considers the time needed to reach a common target. Taking as target the median of the 20 final tour lengths of MA-CETSP, MA-CETSP-BCD reaches it in all 62 instances, in about half the expected running time of MA-CETSP (typical ratio 0.54), and is faster in 47 instances; the 15 instances in which it is slower are mainly those with overlap ratio 0.30. In contrast, V1, although its runs are shorter, needs on average more time than MA-CETSP to reach the same target and does not reach it in 4 instances.

\begin{table}[ht!]
\TABLE
{Computational effort of the variants of Table~\ref{tab_ma_ablation} (20 runs per instance). \label{tab_ma_effort}}
{\begin{tabular}{lcccccc}
\hline
\multirow{2}{*}{Variant} & $\bar t$ & $t_{\mathrm{tot}}$ & ERT & ERT vs.\ V0 & Typical & Target \\
& (s) & (s) & median (s) & F/E/S & ratio & not reached \\
\hline
\hline
V0 (MA-CETSP) & 45.2 & 60.9 & 15.5 & -- & -- & 0 \\
\rowcolor{gray!15}
V1 & 34.9 & 45.6 & 20.2 & 36/1/21 & 1.15 & 4 \\
V2 & 94.9 & 132.7 & 16.9 & 41/0/19 & 1.03 & 2 \\
\rowcolor{gray!15}
V3 & 95.2 & 146.7 & 10.2 & 46/1/15 & 0.59 & 0 \\
V4 (MA-CETSP-BCD) & 80.6 & 123.7 & 7.6 & 47/0/15 & 0.54 & 0 \\
\hline
\end{tabular}}
{$\bar t$: time to the best solution of each run; $t_{\mathrm{tot}}$: time until the stopping criterion; both averaged over the runs and then over the 62 instances. ERT (expected running time): sum of the times spent by the 20 runs until reaching the target, the median of the 20 final tour lengths of V0, or until stopping, divided by the number of runs that reach it; median over the 58 instances in which all variants reach the target. F/E/S: instances in which the ERT is smaller (faster), equal within 5\%, or larger (slower) than that of V0; typical ratio: geometric mean of the ERT ratio to V0 over these instances; target not reached: instances in which no run reaches the target.}
\end{table}
\par

\section{Conclusions} \label{sec5}

In this work, we presented a specialized Block Coordinate Descent (BCD) method to solve the problem of determining visit points for a fixed visiting sequence in the Close Enough Traveling Salesman Problem (CETSP). By exploiting the geometry and algebraic properties of the single-disk subproblem, we demonstrated that optimizing the position of each visit point can be reduced to minimizing a one-dimensional nonlinear function, solved extremely efficiently using Newton's method. Furthermore, we derived an exact algebraic characterization of the stationary points via the roots of a quartic polynomial, and numerically confirmed that the Newton-based approach attains machine-precision accuracy at a substantially lower computational cost than methods based on computing the roots of the quartic polynomial.

Evaluating sequences via the proposed BCD procedure enabled tour evaluations more than an order of magnitude faster than conventional approaches based on Second-Order Cone Programming (SOCP). This drastic reduction in evaluation time made it possible to develop a heuristic method based on Iterated Local Search (ILS) capable of evaluating a significant number of visiting sequences within a combinatorial local search procedure. Testing such a massive volume of candidate sequences would be computationally infeasible when solving the problem via SOCP using generic continuous optimization solvers due to their high evaluation cost. Furthermore, the BCD procedure guarantees disk constraint feasibility at machine precision and naturally exploits warm-starting between neighboring sequences in the search space. Numerical experiments on benchmark instances demonstrated that combining the ILS metaheuristic with the BCD evaluator matched the known optimal values of all small-scale instances and exhibited competitive performance on large-scale problems. Within the memetic algorithm MA-CETSP, one of the best-performing methods for the CETSP, BCD replaced the general-purpose SOCP solver without loss of quality and with a shorter running time, and its use in the local search and after the LKH step led to a method whose best solution over 20 runs is better than that of the original algorithm in 20 of the 62 large instances and equal in the others, that improves the best values published for the original algorithm in 21 instances, and that reaches the solution quality of the original algorithm in about half the time.

As directions for future research, we intend to extend the BCD method to optimize launch and recovery points within the Carrier-Vehicle Traveling Salesman Problem, where the departure and rendezvous locations of the secondary vehicle inside the carrier's regions must be determined jointly. We also plan to investigate adapting the BCD decomposition to CETSP instances featuring non-circular visiting regions, such as ellipses. Finally, integrating the BCD evaluation module as a core building block into other exact or metaheuristic frameworks for routing problems with neighborhoods represents another promising avenue for future work.

\FloatBarrier
\section*{Acknowledgments}
We sincerely thank Marcelo Marazzi for his careful editing and helpful suggestions that enhanced the quality of the written text in English. We also thank Walton P. Coutinho and William K. Mennell for making available the instances used in~\citet{Coutinho2016} and~\citet{mennell2009thesis}, respectively, and Zhenyu Lei and Jin-Kao Hao for making the source code of MA-CETSP~\citep{Lei2024} publicly available.

\bibliographystyle{plainnat}
\bibliography{berma2026}

\FloatBarrier
\begin{appendices}
\section{Algebraic Formulation: Reduction to a Quartic Polynomial} \label{algeform}

This appendix derives the exact algebraic characterization of the stationary points of problem~\eqref{single-disk-1D} summarized in Section~\ref{algesum}. Two different ways are described below.

\subsection{The tangent half-angle substitution}

Since $\hat y_i'(\theta)^T(\hat y_i(\theta)-\hat y_j)$ in~\eqref{firstder} is, geometrically, the two-dimensional cross product $\hat y_j\times\hat y_i(\theta) = (\hat y_j)_x(\hat y_i(\theta))_y - (\hat y_j)_y(\hat y_i(\theta))_x$, multiplying $f'(\theta)=0$ by the product of both norms $\|\hat y_i(\theta)-\hat y_{i-1}\|\,\|\hat y_i(\theta)-\hat y_{i+1}\|$ turns the stationarity condition into the equivalent relation
\begin{equation}
\|\hat y_i(\theta)-\hat y_{i+1}\|\big(\hat y_{i-1}\times\hat y_i(\theta)\big) + \|\hat y_i(\theta)-\hat y_{i-1}\|\big(\hat y_{i+1}\times\hat y_i(\theta)\big) = 0. \label{eq:reflection-relation}
\end{equation}
While the cross-product terms in~\eqref{eq:reflection-relation} are free of square roots, the Euclidean norms $\|\hat y_i(\theta)-\hat y_{i-1}\|$ and $\|\hat y_i(\theta)-\hat y_{i+1}\|$ still involve trigonometric functions inside radicals. Applying the tangent half-angle substitution $t=\tan(\theta/2)$, which expresses $\cos\theta = (1-t^2)/(1+t^2)$ and $\sin\theta = 2t/(1+t^2)$, both radicals become rational functions of $t$ once squared, so squaring~\eqref{eq:reflection-relation} once and clearing denominators already suffices to eliminate them both simultaneously; this produces a degree-$8$ polynomial in $t$ that, with the help of the computer algebra software SymPy~\citep{sympy}, factors as $(1+t^2)\,Q_{\text{spur}}(t)\,Q(t) = 0$, where $1+t^2$ (degree $2$) has no real roots, $Q_{\text{spur}}(t)$ (degree $2$) is an extraneous factor introduced by the squaring step, and $Q(t)$ (degree $4$) is the quartic whose real roots include the stationary points of~\eqref{single-disk-1D}; as discussed below, only two of its four roots are stationary points of $f$. Writing $\hat y_{i-1}=(a_1,b_1)$ and $\hat y_{i+1}=(a_2,b_2)$, we have that
\[
Q(t) = At^4+Bt^3+Ct^2+Dt+E,
\]
with
\[
\begin{aligned}
A &= r_i(b_1+b_2)+(a_1b_2+a_2b_1), \\
B &= 2r_i(a_1+a_2)+4(a_1a_2-b_1b_2), \\
C &= -6(a_1b_2+a_2b_1), \\
D &= 2r_i(a_1+a_2)-4(a_1a_2-b_1b_2), \\
E &= -r_i(b_1+b_2)+(a_1b_2+a_2b_1).
\end{aligned}
\]

\subsection{Reflection principle of geometric optics}

Inspired by the treatment of the Ptolemy--Alhazen mirror-reflection problem in~\citet{fujimura2018}, we represent points as complex numbers instead. Writing $p=\hat y_{i-1}$ and $q=\hat y_{i+1}$ (identifying $\mathbb{R}^2$ with $\mathbb{C}$ in the usual way) and $z=\hat y_i(\theta)=r_ie^{i\theta}$, we characterize the stationary points of $f$ through the classical reflection principle of geometric optics. This equivalence between $f'(\theta)=0$ and the law of reflection can be seen directly from equation~\eqref{firstder} once its common factor $\hat y_i'(\theta)$ is pulled out. Writing $\hat e_1(\theta) = (\hat y_i(\theta)-\hat y_{i-1})/\|\hat y_i(\theta)-\hat y_{i-1}\|$ and $\hat e_2(\theta) = (\hat y_i(\theta)-\hat y_{i+1})/\|\hat y_i(\theta)-\hat y_{i+1}\|$, equation~\eqref{firstder} reads $f'(\theta) = \hat y_i'(\theta)^T[\hat e_1(\theta)+\hat e_2(\theta)]$. Since $\hat y_i'(\theta)$ is tangent to the circle, $f'(\theta)=0$ states that $\hat e_1(\theta)+\hat e_2(\theta)$ is orthogonal to the tangent, hence parallel to the normal at $\hat y_i(\theta)$; and since the sum of two unit vectors is always parallel to the bisector of the angle between them, this is precisely the statement that the normal bisects the angle between the directions to $p$ and to $q$, i.e., the law of reflection.

Fix $z$ on the circle $|z|=r_i$. The tangent line to the circle at $z$ has direction $iz/r_i$, so reflecting a point $w$ about this line gives $z-\frac{z^2}{r_i^2}\overline{(w-z)}$; in particular, the mirror image of $p$ across the tangent at $z$ is
\[
p'(z) = 2z - \frac{z^2\bar p}{r_i^2}.
\]
By the standard unfolding argument of geometric optics, the law of reflection holds at $z$ if and only if $p'(z)$, $z$, and $q$ are collinear, that is,
\[
\operatorname{Im}\!\big[(p'(z)-z)\,\overline{(q-z)}\big] \;=\; 0.
\]
Substituting $p'(z)-z = z - z^2\bar p/r_i^2$, using $\bar z = r_i^2/z$ (valid since $|z|=r_i$) to eliminate $\bar z$, and clearing denominators, this condition reduces to the polynomial equation
\begin{equation}
P(z) \;=\; \alpha z^4 + \beta z^3 + \delta z + \epsilon \;=\; 0, \label{eq:quartic-z}
\end{equation}
with
\[
\alpha = \bar p\,\bar q, \qquad \beta = -r_i^2(\bar p + \bar q), \qquad \delta = r_i^4(p+q), \qquad \epsilon = -r_i^4\,pq,
\]
and no $z^2$ term. Every root of $P$ satisfies $|z|=r_i$ automatically, so every root corresponds to a genuine angle $\theta$, the direction of $z$.

\subsection{Coinciding roots, and the meaning of the spurious pair}

Despite following different routes, $Q(t)$ and $P(z)$ have exactly the same roots. Since $t=\tan(\theta/2)$ and $z=r_ie^{i\theta}$ correspond to the same angle $\theta$, they are related by the classical identity $e^{i\theta}=(1+it)/(1-it)$, i.e., $z=r_i(1+it)/(1-it)$. Substituting this into~\eqref{eq:quartic-z} and clearing denominators gives, after simplification,
\[
P(z)\Big|_{z=r_i(1+it)/(1-it)}\cdot(1-it)^4 \;=\; -2ir_i^4\,Q(t).
\]
Thus $P(z)=0$ and $Q(t)=0$ have exactly the same roots, under the correspondence $t\leftrightarrow\theta\leftrightarrow z$; the two constructions, though motivated differently, characterize the same underlying algebraic curve.

Not every one of these shared roots, however, corresponds to a stationary point of $f$. Exactly two of the four satisfy $f'(\theta)=0$; the other two instead satisfy the equation obtained from~\eqref{firstder} by flipping the sign of its second term---that is, writing $\hat e_1(\theta)$ and $\hat e_2(\theta)$ as above, they satisfy $\hat y_i'(\theta)^T[\hat e_1(\theta)-\hat e_2(\theta)]=0$ and are thus stationary points of $h(\theta):=\|\hat y_i(\theta)-\hat y_{i-1}\|-\|\hat y_i(\theta)-\hat y_{i+1}\|$ rather than of $f$. This is because the collinearity condition above encodes only that $p'(z)$, $z$, and $q$ lie on a common line, and does not distinguish whether $q$ lies on the same ray from $z$ as $p'(z)$ or on the opposite one---the latter being the case that genuinely corresponds to $f'(\theta)=0$, the former instead corresponding to $h'(\theta)=0$. The two cases can nonetheless be distinguished computationally: $\operatorname{Re}[(q-z)/(p'(z)-z)]$ indicates whether $p'(z)$ and $q$ lie on the same side of $z$ (positive) or on opposite sides (negative) whenever the three points are collinear, and we have found, empirically, that the roots $z$ for which this quantity is negative are always exactly the two genuine stationary points of $f$, and the roots for which it is positive are always the stationary points of $h$; we verify this computationally in Section~\ref{evalthree}.

\subsection{Finding the roots}

Every quartic is solvable in radicals via Ferrari's method: writing~\eqref{eq:quartic-z} in depressed form $u^4+pu^2+qu+r=0$ (via a suitable shift in $z$), one solves a resolvent cubic by Cardano's formula and reconstructs the four roots of the quartic from one of its roots via two quadratic solves. The resulting closed-form expression for $\theta^\star$ is, however, impractical: it spans several lines even after simplification, with deeply nested square and cube roots, and is numerically delicate, since Cardano's formula for three real roots requires taking cube roots of complex numbers and discarding the imaginary part, a classic source of catastrophic cancellation if implemented naively. For these reasons, in practice we solve $P(z)$ numerically; the same considerations and methods apply equally to $Q(t)$, given their coinciding roots.

To find the roots of~\eqref{eq:quartic-z} numerically, we use the Durand--Kerner (Weierstrass) method. Given the monic polynomial $P(z)/\alpha = z^4+\tilde\beta z^3+\tilde\delta z+\tilde\epsilon$, with $\tilde\beta=\beta/\alpha$, etc., we choose four initial guesses $z_i^{(0)}$ spread on a circle of radius $R_0=1+\max\{|\tilde\beta|,|\tilde\delta|,|\tilde\epsilon|\}$, which provides a crude bound on the modulus of the roots. The method then iterates
\[
z_i^{(k+1)} = z_i^{(k)} - \frac{P\big(z_i^{(k)}\big)}{\prod_{j\neq i}\big(z_i^{(k)}-z_j^{(k)}\big)}, \qquad i=1,\dots,4,
\]
simultaneously for all four root approximations. Alternatively, we also find the roots via the $4\times 4$ companion matrix
\[
M = \begin{pmatrix}
-\tilde\beta & 0 & -\tilde\delta & -\tilde\epsilon \\
1 & 0 & 0 & 0 \\
0 & 1 & 0 & 0 \\
0 & 0 & 1 & 0
\end{pmatrix},
\]
whose eigenvalues are exactly the roots of $P$. Since $M$ is complex-valued, its eigenvalues are computed using LAPACK's zgeev. (The analogous companion matrix for $Q(t)$ would have $-C/A$, rather than a zero, in its second entry, reflecting the $t^2$ term that $P(z)$ lacks.)

\subsection{Solving $Q(t)$ instead of $P(z)$}

While the roots of $Q(t)$ can vary widely in magnitude and grow unbounded as $\theta$ approaches $\pi$, every root of $P(z)$ lies exactly on the circle $|z|=r_i$. Also, the convergence of the Durand--Kerner method is known to depend strongly on the quality of the initial approximations~\citep{guggenheimer1986}. Therefore, using tighter, more uniform bounds on the roots' magnitude generally yields faster and more stable convergence. For this reason, we solved $P(z)$ rather than $Q(t)$ in the experiments of Section~\ref{evalthree}. Separate numerical experiments applying Durand--Kerner to $Q(t)$ instead of $P(z)$ show that convergence takes on average $61.2$ iterations and $3.93\times10^{-6}$ seconds per call --- respectively $2.8\times$ and $2.5\times$ slower than solving $P(z)$. The companion-matrix approach shows the opposite pattern: applied to $Q(t)$, whose real coefficients allow the real eigensolver dgeev instead of the complex zgeev, it takes $1.69\times10^{-6}$ seconds per call, $3.5\times$ faster than applied to $P(z)$. Regardless of which polynomial and which root-finding method is used, however, the adaptive multi-start Newton scheme remains the fastest by a wide margin: $8.2\times$ faster than Durand--Kerner on $Q(t)$ and $3.5\times$ faster than the companion-matrix approach on $Q(t)$, in addition to the $3.3\times$ and $12.4\times$ margins reported in Section~\ref{secadapt} for $P(z)$.

\section{Analysis of the Anchors of the Multi-Start Newton Method} \label{appadapt}

This appendix details the analysis behind the adaptive multi-start procedure of Section~\ref{secadapt}, carried out on the dataset of Section~\ref{evalthree}.

In the multi-start 1D Newton, the three anchors converge to the global minimizer at nearly identical individual rates --- $\theta_0^1$, $\theta_0^2$, and $\theta_0^3$ hit $\theta^\star$ in $99.49\%$, $99.50\%$, and $99.32\%$ of instances, respectively --- so no single anchor stands out as more reliable on its own. The asymmetry instead shows up in how the anchors' failures overlap. Considering them in pairs, $\theta_0^1$ and $\theta_0^2$ together cover $99.9838\%$ of instances, whereas $\theta_0^1$ with $\theta_0^3$ covers only $99.6853\%$, and $\theta_0^2$ with $\theta_0^3$ only $99.6937\%$. Thus, although $\theta_0^3$ hits $\theta^\star$ about as often as $\theta_0^1$ or $\theta_0^2$ individually, its failures coincide with theirs far more often than $\theta_0^1$'s and $\theta_0^2$'s failures coincide with each other.

We now explain this asymmetry between $\theta_0^3$ and $\{\theta_0^1,\theta_0^2\}$. A first-order Taylor expansion of~$f$ in~$r_i$ around $r_i=0$, with $y(\theta) = c_i + r_i(\cos\theta,\sin\theta)$ and $y_{i-1},y_{i+1}$ held fixed, gives
\begin{equation} \label{taylor}
f(\theta) = \|\hat y_{i-1} \| + \| \hat y_{i+1} \| - r_i \|\hat u_1+\hat u_2\| \cos(\theta-\theta_0^3) + O\!\left(r_i^2\left(\frac{1}{\|\hat y_{i-1} \|}+\frac{1}{\|\hat y_{i+1} \|}\right)\right),
\end{equation}
where the leading term follows from the identity $\cos(\theta-\theta_0^1)+\cos(\theta-\theta_0^2)=\|\hat u_1 + \hat u_2\|\cos(\theta-\theta_0^3)$. The approximation of $f$, disregarding the remainder, attains its minimum at $\theta_0^3$ and its maximum at $\theta_0^3+\pi$. What matters for $\theta_0^3$ to reliably approximate the minimizer $\theta^\star$ of $f$ is the size of the remainder relative to the amplitude of this leading term, $r_i\|\hat u_1+\hat u_2\|$; the two are related by
\[
\frac{\text{remainder}}{\text{leading amplitude}} = O\!\left(\frac{r_i}{\|\hat u_1+\hat u_2\|}\left(\frac{1}{\|\hat y_{i-1} \|}+\frac{1}{\|\hat y_{i+1} \|}\right)\right),
\]
which is small precisely when $r_i$ is small relative to $\|y_{i-1}-c_i\|$ and $\|y_{i+1}-c_i\|$. (This can fail when $\hat u_1\approx-\hat u_2$ (i.e., $y_{i-1}$ and $y_{i+1}$ nearly antipodal about $c_i$), in which case $\|\hat u_1+\hat u_2\|\to0$ and the leading-order term itself vanishes, so the approximation offers no guarantee regardless of how small $r_i$ is. This degenerate case cannot occur, however: $\hat u_1=-\hat u_2$ would place $c_i$ on the open segment $(y_{i-1},y_{i+1})$, in which case the segment intersects the disk.) Consequently, $\theta_0^3$ is essentially guaranteed to reach the global minimizer whenever $f$ is flat. On the other hand, $\theta_0^1$ and $\theta_0^2$ have no privileged relationship to the leading-order term and succeed only through the remainder --- an effect that becomes proportionally smaller, and correspondingly harder for Newton's method to resolve, as $f$ flattens.

This connection between flatness and the reliability of $\theta_0^3$ as a starting point for Newton suggests that, in practice, $\theta_0^3$ should be used only when $f$ is flat --- precisely the regime in which Newton's ability to converge starting from $\theta_0^1$ or $\theta_0^2$ is compromised, and in which $\theta_0^3$ performs well. This can be tested directly: define, for each instance, $\rho_i = r_i\|\hat u_1+\hat u_2\| / (\|\hat y_{i-1}\|+\|\hat y_{i+1}\|)$, the relative peak-to-peak amplitude of $f$ predicted by the first-order expansion~\eqref{taylor}. Splitting the $1{,}000{,}000$ instances by the joint convergence status of $\theta_0^1$ and $\theta_0^2$, the two groups are cleanly separated: in the $990{,}121$ instances where both converge to $\theta^\star$, $\rho_i$ has median $0.269$ and $5$th percentile $0.036$; in the $162$ instances where both fail, $\rho_i$ has median $0.0068$ and never exceeds $0.0085$ --- below the $5$th percentile of the successful group. Every instance in which $\theta_0^1$ and $\theta_0^2$ jointly fail is thus among the flattest in the dataset, confirming that it is precisely the near-flatness of $f$ that makes the line-of-sight anchors unreliable. A threshold as low as $\rho_i < 0.01$ would already separate the two groups correctly in this dataset, but we adopt the more conservative practical rule of using $\theta_0^3$ whenever $\rho_i < 0.05$, a nearly sixfold margin above the largest $\rho_i$ observed among the $162$ failing instances, to guard against instances outside this particular sample.

Independently of flatness, there is also a substantial overlap between the instances in which Newton converges starting from $\theta_0^1$ and those in which it converges starting from $\theta_0^2$, suggesting that starting from both is not always necessary either. Among the $9{,}717$ instances in which exactly one of $\theta_0^1$ and $\theta_0^2$ reaches $\theta^\star$, choosing a starting point using a criterion that does not involve the Newton method at all --- namely, start from whichever of $\theta_0^1$ and $\theta_0^2$ has the smaller value of $f$ --- picks the correct starting point in $9{,}630$ cases ($99.10\%$). The $87$ cases in which this predictor fails are themselves cleanly separated from the rest by the relative gap $|f(\theta_0^1)-f(\theta_0^2)|/\min(f(\theta_0^1),f(\theta_0^2))$: in the $9{,}630$ correctly-predicted instances this gap has median $0.125$ and $5$th percentile $0.027$, whereas in the $87$ mispredicted instances it never exceeds $2.4\times10^{-5}$. We adopt the practical rule of trusting the smaller-$f$ anchor alone whenever this gap exceeds $10^{-4}$, a more than fourfold margin above the largest gap observed among the $87$ mispredicted instances, and trying both when it does not. Combining the two rules yields the adaptive multi-start procedure described in Section~\ref{secadapt}.

\section{Neighborhood Structures of the VND} \label{appvnd}

The VND of Section~\ref{sec43} considers eight neighborhood structures. The 2-opt and 3-opt neighborhoods remove two and three non-consecutive tour edges, respectively, and reconnect the resulting segments to form another tour. The 3-opt neighborhood excludes moves already covered by 2-opt. The double-bridge neighborhood uses structured 4-opt moves and is explored through Glover's approach~\citep{Glover1996}, which models the search over a restricted class of reconnections as a shortest-path problem in a directed acyclic graph. The swap neighborhood exchanges two non-adjacent customers, whereas the insertion neighborhood relocates a single customer to another position in the tour. The string relocation neighborhood moves a segment of consecutive customers to another position, optionally reversing its order upon reinsertion. Segments of lengths from two to a maximum length $L$ are considered. The Balas--Simonetti neighborhood~\citep{Balas2001} depends on a parameter $k$: for positions $a,b\in\{1,\dots,n\}$ satisfying $b\geq a+k$, customer $i_a$ must precede customer $i_b$ in the new sequence. The resulting restricted reordering problem is explored using an auxiliary directed acyclic graph. The Gutin neighborhood~\citep{Gutin2005} selects a set $Z\subseteq\{1,\dots,n\}$ of pairwise non-adjacent tour positions, with adjacency understood cyclically, and reassigns the customers occupying those positions. With the visit points of the remaining customers held fixed, this reassignment is formulated as an assignment problem. For $a,b\in Z$, the cost of assigning customer $i_a$ to position $b$ is
$$
C_{ab} = \min_{u\in B(c_{i_a},r_{i_a})} \left\{ \|u-y_{i_{b-1}}\|_2 + \|u-y_{i_{b+1}}\|_2 \right\}.
$$
Since the selected positions are pairwise non-adjacent, the visit points immediately before and after each target position remain fixed when computing these costs. Each entry $C_{ab}$ is therefore obtained by solving a single-disk subproblem. After solving the assignment problem, BCD recomputes all visit points for the resulting sequence. The parameter $|Z|$ controls the number of reassigned customers. As in~\citet{ber2025-tspn}, the 2-opt, 3-opt, swap, insertion, and string relocation neighborhoods are explored with candidate lists: for each customer we keep the $K$ disks closest to it, and only the moves that insert an edge between a customer and one of its candidates are considered. The lists are built once, before the search, from the distance between disks, $\max\{0,\|c_i-c_j\|_2-r_i-r_j\}$. The parameter $K$ controls their size.

\section{Calibration of the ILS Parameters} \label{seccal}

The ILS metaheuristic described in Section~\ref{sec43} has thirteen parameters to be calibrated. Three govern acceptance: the initial threshold $\eta_0$, the reduction factor $h$, and the number of reductions $\kappa$. Six concern VND: whether the neighborhood exploration order is fixed or randomly chosen at each call; whether 2-opt, 3-opt, swap, and insertion use first improvement or best improvement; the Balas--Simonetti parameter $k$; the number $|Z|$ of positions reassigned by the Gutin neighborhood; the maximum string length $L$; and whether relocated strings may be reversed. Two concern the computation of visit points within VND: the maximum number $s_{\max}$ of BCD sweeps performed when evaluating each move, and whether a complete BCD (up to 200 sweeps) is applied to the solution returned by VND. The remaining two parameters are the perturbation strength $\ell$ and the candidate-list size $K$, which restricts the searches in the 2-opt, 3-opt, swap, insertion, and string relocation neighborhoods.

All thirteen parameters were calibrated simultaneously using \texttt{Irace}~\citep{LopezIbanez2016} (Iterated Racing, version~4.3), with elitist racing and a budget of 1,500 runs. The remaining racing settings were left at their defaults, including a confidence level of 0.95 and the settings automatically determined from the budget and the number of parameters being calibrated, such as the number of racing iterations and the minimum number of surviving configurations. Each algorithm run during calibration was subject to an external time limit of 600~seconds.

The training set consisted of 20 instances selected from the CETSP instances generated in Section~\ref{sec42}, with numbers of customers $n\in\{10,40,70,90\}\cup\{100,120,140,\dots,400\}$. One instance per size was selected to avoid biasing the calibration toward a particular instance size. An independent test set consisting of 20 additional instances from the same collection, with the same sizes but different disk configurations, was used to assess the elite configurations on instances outside the training set. The calibration procedure terminated after six racing iterations, when the remaining budget (11 of the 1,500 runs) was no longer sufficient for a new race, and identified five elite configurations. The four best of them were evaluated on the test set, and their performance was similar: for each test instance, the tour length obtained by any of these configurations differed by less than $2.12\%$ from the best tour length obtained by the four configurations.

All five elite configurations used first improvement, $s_{\max}=1$, and a complete BCD applied to the solution returned by VND. An analysis of the test results by instance size showed that the configuration ranked first by \texttt{Irace} also achieved the best aggregate performance on the large instances ($n \geq 260$ disks), yielding the smallest sum of tour lengths in this range, as well as the smallest average tour length over the whole test set. Since performance on large instances is the primary criterion for the experiments in Section~\ref{sec43}, this configuration was adopted for all subsequent experiments. Its parameter values are $\kappa=20$, $\eta_0=0.02$, $h=0.85$, $L=3$, $|Z|=10$, $\ell=5$, $k=3$, $K=16$ and $s_{\max}=1$, with string reversal disabled, first improvement used in the 2-opt, 3-opt, swap, and insertion neighborhoods, the neighborhoods explored in the fixed order in which they were listed, rather than in an order drawn at random at each call of the VND, and a complete BCD applied to the solution returned by VND.

\section{Detailed Per-Instance Results} \label{appresults}
\setcounter{table}{0}
\renewcommand{\thetable}{\thesection.\arabic{table}}

This appendix reports the instance-by-instance results summarized in Section~\ref{secmeta}: the tour lengths and computational effort of ILS-BCD (Tables~\ref{tab_ils_effort} and~\ref{tab_ils_effort_var}), the tour lengths obtained by ILS-BCD and by the heuristics evaluated in~\citet{mennell2009thesis} (Tables~\ref{tab_cetsp2.1} and~\ref{tab_cetsp4.1}), by four more recent heuristics (Tables~\ref{tab_recent_common} and~\ref{tab_recent_var}), and by MA-CETSP and MA-CETSP-BCD (Table~\ref{tab_ma_instances}).

\begin{table}[ht!]
\TABLE
{Results of ILS-BCD on the 48 benchmark instances of~\citet{mennell2009thesis} with a common radius: tour length and computational effort until the best solution was found. \label{tab_ils_effort}}
{\resizebox{0.60\fullwidth}{!}{
\begin{tabular}{ccccccc}
\hline
& & & & \multicolumn{3}{c}{Computational effort} \\
\cmidrule(lr){5-7}
Instance & $n$ & OR & Tour Length & CPU Time (s) & ILS iter & BCD calls \\
\hline
\hline
           kroD100 &   100 & 0.02 &   159.037 &    343.50 &    767 &    11,303,663 \\
\rowcolor{gray!15}
            rat195 &   195 & 0.02 &   157.654 &  7,862.88 &  2,256 &   134,757,937 \\
            lin318 &   318 & 0.02 & 2,830.558 & 29,115.53 & 11,195 &   391,912,101 \\
\rowcolor{gray!15}
             rd400 &   400 & 0.02 & 1,026.187 & 27,231.78 &  7,427 &   247,908,940 \\
            pcb442 &   442 & 0.02 &   320.658 & 27,272.98 &  5,299 &   294,271,966 \\
\rowcolor{gray!15}
              d493 &   493 & 0.02 &   202.455 & 35,355.06 &  2,866 &   382,201,820 \\
           dsj1000 & 1,000 & 0.02 &   969.346 & 35,962.87 &  1,191 &   188,452,815 \\
\rowcolor{gray!15}
           kroD100 &   100 & 0.10 &    89.668 &  3,369.18 &  9,972 &   226,502,938 \\
            rat195 &   195 & 0.10 &    67.991 & 15,775.75 &  7,827 &   628,447,072 \\
\rowcolor{gray!15}
            lin318 &   318 & 0.10 & 1,394.627 & 13,460.90 &  2,539 &   361,177,915 \\
             rd400 &   400 & 0.10 &   466.082 & 33,084.78 &  4,906 &   654,507,123 \\
\rowcolor{gray!15}
            pcb442 &   442 & 0.10 &   142.716 & 32,805.76 &  3,039 &   664,817,122 \\
              d493 &   493 & 0.10 &   100.802 & 33,346.40 &    811 &   617,371,399 \\
\rowcolor{gray!15}
           dsj1000 & 1,000 & 0.10 &   377.617 & 25,950.93 &    567 &   204,982,259 \\
           kroD100 &   100 & 0.30 &    58.541 &  1,676.56 &  2,923 &   155,302,529 \\
\rowcolor{gray!15}
            rat195 &   195 & 0.30 &    45.702 &  1,165.30 &      9 &    77,191,009 \\
            lin318 &   318 & 0.30 &   765.980 & 33,465.02 &     97 & 1,359,488,134 \\
\rowcolor{gray!15}
             rd400 &   400 & 0.30 &   224.841 & 33,234.49 &  1,555 &   809,923,865 \\
            pcb442 &   442 & 0.30 &    83.540 & 17,889.29 &    942 &   309,256,149 \\
\rowcolor{gray!15}
              d493 &   493 & 0.30 &    69.759 & 31,756.01 &     67 &   849,042,574 \\
           dsj1000 & 1,000 & 0.30 &   199.972 & 15,168.85 &     89 &    91,022,897 \\
\rowcolor{gray!15}
         team1-100 &   101 &   -- &   307.337 &  2,065.17 &  7,423 &   119,035,700 \\
         team2-200 &   201 &   -- &   246.683 &  5,436.14 &  2,261 &   234,349,460 \\
\rowcolor{gray!15}
         team3-300 &   301 &   -- &   461.889 & 34,222.41 &  9,572 &   769,333,048 \\
         team4-400 &   401 &   -- &   680.090 & 31,973.98 &  5,686 &   449,787,890 \\
\rowcolor{gray!15}
         team5-499 &   500 &   -- &   708.850 & 35,634.23 &  3,595 &   370,834,164 \\
         team6-500 &   501 &   -- &   225.241 & 35,023.28 &  1,341 &   604,479,096 \\
\rowcolor{gray!15}
         bonus1000 & 1,001 &   -- &   397.188 & 35,404.47 &    797 &   268,784,268 \\
 rotatingDiamonds1 &    21 &   -- &    32.389 &      0.00 &      0 &            21 \\
\rowcolor{gray!15}
 rotatingDiamonds2 &    61 &   -- &   140.477 &     37.62 &    396 &     3,936,276 \\
 rotatingDiamonds3 &   181 &   -- &   380.882 &    797.76 &  1,221 &    37,642,380 \\
\rowcolor{gray!15}
 rotatingDiamonds4 &   321 &   -- &   770.661 & 16,230.05 &  8,938 &   468,336,159 \\
 rotatingDiamonds5 &   681 &   -- & 1,510.753 & 29,772.26 &  6,046 &   437,839,110 \\
\rowcolor{gray!15}
          bubbles1 &    37 &   -- &   349.135 &      0.34 &     12 &        53,485 \\
          bubbles2 &    77 &   -- &   428.279 &    136.12 &    809 &    14,767,122 \\
\rowcolor{gray!15}
          bubbles3 &   127 &   -- &   529.955 &  1,218.63 &  1,510 &    43,123,709 \\
          bubbles4 &   185 &   -- &   803.393 &  7,439.88 &  4,488 &   195,584,273 \\
\rowcolor{gray!15}
          bubbles5 &   251 &   -- & 1,035.825 & 28,036.78 &  9,310 &   547,002,355 \\
          bubbles6 &   325 &   -- & 1,227.248 & 21,305.28 &  5,324 &   334,908,363 \\
\rowcolor{gray!15}
          bubbles7 &   407 &   -- & 1,578.313 & 31,836.17 &  4,271 &   350,701,644 \\
          bubbles8 &   497 &   -- & 1,938.889 & 35,546.95 &  3,476 &   323,682,071 \\
\rowcolor{gray!15}
          bubbles9 &   595 &   -- & 2,317.461 & 35,410.64 &  2,312 &   259,367,372 \\
concentricCircles1 &    17 &   -- &    53.158 &      0.00 &      0 &            17 \\
\rowcolor{gray!15}
concentricCircles2 &    37 &   -- &   153.132 &      1.92 &     21 &       142,928 \\
concentricCircles3 &    61 &   -- &   270.007 &    353.42 &  1,658 &    18,141,533 \\
\rowcolor{gray!15}
concentricCircles4 &   105 &   -- &   451.870 &  1,334.98 &  1,777 &    39,673,013 \\
concentricCircles5 &   149 &   -- &   632.978 & 15,284.15 & 16,152 &   341,823,306 \\
\rowcolor{gray!15}
     chaoSingleDep &   201 &   -- & 1,039.611 & 26,708.41 & 84,776 &   950,144,983 \\
\hline
           Average &   337 &   -- &   592.196 & 18,468.85 &  5,198 &   330,069,083 \\
\hline
\end{tabular}}}
{}
\end{table}

\begin{table}[ht!]
\TABLE
{Results of ILS-BCD on the 14 benchmark instances of~\citet{mennell2009thesis} with variable radii: tour length and computational effort until the best solution was found. \label{tab_ils_effort_var}}
{\resizebox{0.60\fullwidth}{!}{
\begin{tabular}{cccccc}
\hline
& & & \multicolumn{3}{c}{Computational effort} \\
\cmidrule(lr){4-6}
Instance & $n$ & Tour Length & CPU Time (s) & ILS iter & BCD calls \\
\hline
\hline
  kroD100 &   100 &   141.829 &    120.59 &    502 &     4,783,382 \\
\rowcolor{gray!15}
   rat195 &   195 &    68.224 & 28,389.70 &  1,228 & 1,526,290,647 \\
   lin318 &   318 & 2,078.281 &  3,298.02 &    653 &    68,379,654 \\
\rowcolor{gray!15}
    rd400 &   400 & 1,243.392 & 34,811.18 &  4,552 &   213,711,425 \\
   pcb442 &   442 &   219.842 & 34,006.94 &  4,581 &   650,581,586 \\
\rowcolor{gray!15}
     d493 &   493 &   134.258 & 22,786.57 &  1,177 &   362,010,265 \\
  dsj1000 & 1,000 &   636.998 & 32,285.83 &    841 &   248,778,878 \\
\rowcolor{gray!15}
team1-100 &   101 &   388.537 &    162.81 &    613 &     7,993,587 \\
team2-200 &   201 &   613.659 & 13,080.35 & 11,382 &   277,385,488 \\
\rowcolor{gray!15}
team3-300 &   301 &   378.087 & 34,211.33 &  7,935 &   925,696,285 \\
team4-400 &   401 & 1,000.134 & 35,961.59 &  5,330 &   314,266,843 \\
\rowcolor{gray!15}
team5-499 &   500 &   446.192 & 32,968.15 &  1,490 &   611,726,108 \\
team6-500 &   501 &   640.234 & 34,883.96 &  2,512 &   493,706,471 \\
\rowcolor{gray!15}
bonus1000 & 1,001 &   980.031 & 35,285.60 &    731 &   209,950,717 \\
\hline
  Average &   425 &   640.693 & 24,446.62 &  3,109 &   422,518,667 \\
\hline
\end{tabular}}}
{}
\end{table}

\begin{table}[ht!]
\TABLE
{Tour lengths obtained by ILS-BCD and by the 15 heuristics evaluated in~\citet{mennell2009thesis} on the 48 benchmark instances with a common radius. \label{tab_cetsp2.1}}
{\resizebox{\fullwidth}{!}{
\begin{tabular}{cccccccccccccccccc}
\hline
& & \multicolumn{15}{c}{Mennell (2009), Table 2.2 (pp. 37–38)} & \\
\cmidrule(lr){3-17}
Instance & ILS-BCD & HYBRID2 & HYBRID1 & SZ2 & SZ1 & CTD & GTSP2 & YUAN+ & LK-SOCP & TILE3 & TILE4 & TILE2 & TILE1 & GTSP1 & MERGE & HULL & Best \\
\hline
\hline
             kroD100 &   \textbf{159.04} &         160.46 &   160.88 &         163.76 &         163.78 &   165.34 &          159.05 &            164.31 &            164.31 &   170.88 &   170.32 &   172.84 &   173.05 &   159.28 &   188.07 &   210.03 &   159.04 \\
\rowcolor{gray!15}
              rat195 &   \textbf{157.65} &         164.51 &   165.37 &         171.29 &         171.35 &   169.33 &          158.79 &            170.89 &            169.32 &   181.99 &   188.76 &   187.37 &   190.00 &   159.09 &   189.12 &   229.23 &   157.65 \\
              lin318 & \textbf{2,830.56} &       2,903.26 & 2,914.47 &       2,941.00 &       2,941.53 & 2,968.92 &        2,867.46 &          2,949.71 &          3,031.48 & 3,191.60 & 3,155.60 & 3,217.70 & 3,183.20 & 2,875.34 & 3,261.49 & 3,824.50 & 2,830.56 \\
\rowcolor{gray!15}
               rd400 & \textbf{1,026.19} &       1,039.27 & 1,042.62 &       1,066.16 &       1,066.44 & 1,074.39 &        1,033.42 &          1,069.26 &          1,091.84 & 1,164.00 & 1,166.60 & 1,180.50 & 1,177.90 & 1,035.95 & 1,243.61 & 1,513.97 & 1,026.19 \\
              pcb442 &   \textbf{320.66} &         327.26 &   328.72 &         333.24 &         333.27 &   332.93 &          323.03 &            340.59 &            348.15 &   359.93 &   366.64 &   371.53 &   359.47 &   324.21 &   401.09 &   450.22 &   320.66 \\
\rowcolor{gray!15}
                d493 &   \textbf{202.45} &         209.82 &   211.05 &         213.51 &         213.64 &   214.16 &          204.71 &            221.41 &            217.03 &   227.55 &   236.76 &   230.83 &   231.85 &   206.77 &   249.62 &   261.70 &   202.45 \\
             dsj1000 &            969.35 &         998.35 & 1,015.30 &       1,010.86 &       1,012.72 & 1,002.42 & \textbf{935.74} &          1,005.80 &          1,023.11 & 1,156.20 & 1,140.00 & 1,156.60 & 1,141.10 &   951.09 & 1,266.38 & 1,456.97 &   935.74 \\
\rowcolor{gray!15}
             kroD100 &    \textbf{89.67} &          96.68 &    97.58 &          96.68 &          97.09 &    96.20 &           89.94 &             95.64 &             95.64 &   116.63 &   114.81 &   102.49 &   115.81 &    91.53 &   121.88 &   107.85 &    89.67 \\
              rat195 &    \textbf{67.99} &          73.89 &    74.41 &          74.16 &          74.22 &    76.67 &           68.26 &             77.66 &             82.36 &    82.34 &    84.45 &    87.74 &    84.12 &    72.25 &    90.81 &    85.83 &    67.99 \\
\rowcolor{gray!15}
              lin318 & \textbf{1,394.63} &       1,554.59 & 1,563.48 &       1,556.82 &       1,557.47 & 1,526.39 &        1,467.02 &          1,557.07 &          1,532.05 & 1,842.00 & 1,679.50 & 1,711.80 & 1,646.80 & 1,555.56 & 1,828.95 & 1,886.08 & 1,394.63 \\
               rd400 &   \textbf{466.08} &         521.91 &   526.96 &         521.91 &         526.73 &   521.67 &          473.70 &            528.13 &            541.22 &   602.22 &   601.34 &   584.86 &   623.32 &   498.51 &   608.71 &   593.45 &   466.08 \\
\rowcolor{gray!15}
              pcb442 &   \textbf{142.72} &         161.77 &   162.61 &         162.26 &         162.73 &   163.07 &          147.24 &            214.50 &            211.66 &   179.91 &   186.02 &   185.63 &   180.30 &   159.06 &   197.41 &   177.65 &   142.72 \\
                d493 &   \textbf{100.80} &         110.75 &   110.98 &         110.90 &         111.03 &   108.22 &          112.55 &            114.16 &            116.98 &   119.74 &   115.65 &   115.05 &   121.24 &   133.81 &   124.89 &   112.37 &   100.80 \\
\rowcolor{gray!15}
             dsj1000 &   \textbf{377.62} &         404.60 &   408.69 &         404.60 &         407.02 &   404.89 &          482.85 &            464.26 &            420.33 &   429.90 &   445.21 &   455.71 &   447.28 &   632.79 &   515.66 &   462.72 &   377.62 \\
             kroD100 &    \textbf{58.54} & \textbf{58.54} &    58.69 & \textbf{58.54} & \textbf{58.54} &    58.70 &           62.15 &    \textbf{58.54} &    \textbf{58.54} &    61.45 &    61.45 &    61.45 &    61.45 &    71.32 &    78.45 &    64.55 &    58.54 \\
\rowcolor{gray!15}
              rat195 &    \textbf{45.70} &          45.79 &    45.91 &          45.79 &          45.91 &    45.82 &           48.19 &    \textbf{45.70} &             45.74 &    52.79 &    52.97 &    52.11 &    51.82 &    75.04 &    58.15 &    52.92 &    45.70 \\
              lin318 &   \textbf{765.98} &         772.61 &   780.56 &         782.40 &         782.96 &   803.37 &          946.67 &            776.58 &            776.58 &   809.06 &   809.06 &   809.06 &   809.06 & 1,304.03 & 1,116.14 &   814.65 &   765.98 \\
\rowcolor{gray!15}
               rd400 &   \textbf{224.84} &         233.29 &   237.59 &         233.29 &         233.75 &   239.81 &          348.99 &            228.95 &            225.63 &   300.32 &   300.32 &   294.05 &   298.04 &   445.58 &   255.27 &   239.87 &   224.84 \\
              pcb442 &    \textbf{83.54} &          85.12 &    85.18 &          85.12 &          85.18 &    85.61 &          126.21 &            100.72 &            106.32 &    93.64 &    91.71 &    93.59 &    93.59 &   155.40 &    89.82 &    88.77 &    83.54 \\
\rowcolor{gray!15}
                d493 &    \textbf{69.76} &          70.75 &    70.84 &          70.75 &          70.84 &    75.38 &           82.83 &    \textbf{69.76} &    \textbf{69.76} &    80.09 &    82.89 &    80.71 &    78.94 &   148.97 &    76.95 &    73.19 &    69.76 \\
             dsj1000 &            199.97 &         201.92 &   203.08 &         201.92 &         201.95 &   207.16 &          459.22 &            204.16 &   \textbf{199.95} &   214.38 &   214.38 &   214.87 &   214.87 &   754.84 &   277.26 &   219.71 &   199.95 \\
\rowcolor{gray!15}
           team1-100 &   \textbf{307.34} &         317.09 &   320.10 &         326.73 &         327.28 &   325.89 &          307.79 &            316.32 &            316.32 &   352.34 &   352.34 &   350.27 &   350.27 &   311.90 &   391.71 &   436.15 &   307.34 \\
           team2-200 &   \textbf{246.68} &         277.61 &   283.50 &         277.61 &         277.75 &   287.23 &          249.92 &            278.24 &            345.41 &   305.46 &   304.36 &   325.74 &   305.46 &   286.81 &   322.94 &   334.91 &   246.68 \\
\rowcolor{gray!15}
           team3-300 &   \textbf{461.89} &         504.30 &   508.47 &         505.65 &         508.75 &   503.87 &          484.70 &            497.66 &            497.66 &   576.73 &   572.13 &   557.06 &   558.30 &   501.02 &   586.83 &   610.66 &   461.89 \\
           team4-400 &   \textbf{680.09} &         739.28 &   745.42 &         742.83 &         745.58 &   734.16 &          680.21 &            743.28 &            777.74 &   919.34 &   877.68 &   851.85 &   898.53 &   689.95 &   869.22 &   974.01 &   680.09 \\
\rowcolor{gray!15}
           team5-499 &            708.85 &         727.13 &   730.07 &         727.57 &         727.84 &   729.52 & \textbf{703.20} &            740.67 &            758.48 &   787.30 &   803.12 &   814.46 &   808.53 &   707.47 &   892.65 &   944.55 &   703.20 \\
           team6-500 &   \textbf{225.24} &         232.94 &   235.91 &         232.94 &         233.02 &   253.27 &          355.02 &            282.76 &            245.17 &   266.33 &   266.33 &   266.84 &   266.49 &   487.73 &   278.65 &   260.07 &   225.24 \\
\rowcolor{gray!15}
           bonus1000 &   \textbf{397.19} &         435.98 &   441.00 &         439.59 &         444.43 &   436.42 &          578.87 &            485.30 &            508.27 &   504.98 &   485.03 &   495.37 &   524.21 &   699.33 &   557.04 &   470.73 &   397.19 \\
   rotatingDiamonds1 &    \textbf{32.39} &          33.15 &    33.15 &          33.15 &          33.15 &    33.15 &  \textbf{32.39} &    \textbf{32.39} &    \textbf{32.39} &    32.56 &    32.77 &    34.30 &    35.62 &    32.41 &    38.00 &    34.36 &    32.39 \\
\rowcolor{gray!15}
   rotatingDiamonds2 &   \textbf{140.48} &         141.80 &   141.81 &         142.25 &         142.25 &   142.01 & \textbf{140.48} &            144.71 &            144.71 &   147.09 &   145.39 &   145.51 &   145.51 &   140.53 &   155.24 &   150.76 &   140.48 \\
   rotatingDiamonds3 &   \textbf{380.88} &         382.37 &   382.56 &         383.35 &         383.38 &   386.27 &          382.17 &            385.53 &            385.53 &   386.16 &   386.59 &   386.22 &   385.29 &   383.97 &   420.59 &   430.78 &   380.88 \\
\rowcolor{gray!15}
   rotatingDiamonds4 &   \textbf{770.66} &         772.83 &   773.24 &         774.80 &         774.87 &   778.29 &          773.21 &            771.31 &            771.31 &   779.85 &   776.42 &   777.12 &   773.91 &   783.49 &   814.39 &   882.19 &   770.66 \\
   rotatingDiamonds5 & \textbf{1,510.75} &       1,512.98 & 1,516.65 &       1,515.44 &       1,515.59 & 1,520.82 &        1,517.70 &          1,512.48 &          1,512.49 & 1,515.00 & 1,513.30 & 1,516.30 & 1,525.26 & 1,554.60 & 1,762.00 & 1,822.56 & 1,510.75 \\
\rowcolor{gray!15}
            bubbles1 &   \textbf{349.13} &         360.87 &   360.89 &         351.58 &         351.58 &   361.97 & \textbf{349.13} &   \textbf{349.13} &   \textbf{349.13} &   367.33 &   367.33 &   367.33 &   368.96 &   349.55 &   382.64 &   376.81 &   349.13 \\
            bubbles2 &   \textbf{428.28} &         428.84 &   428.87 &         429.29 &         428.96 &   430.28 &          432.38 &            514.14 &            514.14 &   470.38 &   467.02 &   455.90 &   467.68 &   442.32 &   501.75 &   459.18 &   428.28 \\
\rowcolor{gray!15}
            bubbles3 &   \textbf{529.95} &         590.39 &   591.31 &         595.56 &         595.96 &   565.97 &          546.73 &            712.74 &            712.74 &   609.26 &   614.91 &   682.58 &   725.77 &   565.82 &   814.60 &   933.82 &   529.95 \\
            bubbles4 &   \textbf{803.39} &         843.41 &   844.11 &         863.49 &         863.56 &   915.86 &          877.47 &          1,024.51 &          1,024.51 &   984.69 &   992.59 & 1,025.70 &   971.89 &   901.36 & 1,157.24 & 1,253.66 &   803.39 \\
\rowcolor{gray!15}
            bubbles5 & \textbf{1,035.82} &       1,120.50 & 1,121.71 &       1,141.87 &       1,142.85 & 1,174.62 &        1,131.09 &          1,327.46 &          1,327.46 & 1,278.80 & 1,296.50 & 1,291.70 & 1,302.30 & 1,169.65 & 1,558.71 & 1,479.07 & 1,035.82 \\
            bubbles6 & \textbf{1,227.25} &       1,444.69 & 1,447.05 &       1,464.27 &       1,464.47 & 1,508.27 &        1,446.26 &          1,674.06 &          1,674.06 & 1,611.00 & 1,659.70 & 1,656.70 & 1,579.00 & 1,486.76 & 1,962.85 & 1,760.09 & 1,227.25 \\
\rowcolor{gray!15}
            bubbles7 & \textbf{1,578.31} &       1,735.89 & 1,737.95 &       1,839.94 &       1,839.70 & 1,841.19 &        1,775.92 &          2,154.57 &          2,154.57 & 2,150.60 & 2,182.40 & 2,224.40 & 2,152.80 & 1,878.28 & 2,537.17 & 2,356.86 & 1,578.31 \\
            bubbles8 & \textbf{1,938.89} &       2,136.39 & 2,138.06 &       2,250.56 &       2,251.56 & 2,244.85 &        2,236.40 &          2,630.22 &          2,630.22 & 2,535.50 & 2,476.00 & 2,553.80 & 2,538.40 & 2,353.03 & 2,894.47 & 2,729.28 & 1,938.89 \\
\rowcolor{gray!15}
            bubbles9 & \textbf{2,317.46} &       2,547.50 & 2,537.80 &       2,626.10 &       2,648.08 & 2,684.92 &        2,748.93 &          3,176.34 &          3,176.34 & 2,904.20 & 2,893.90 & 2,850.20 & 2,943.90 & 2,988.66 & 3,727.61 & 3,368.31 & 2,317.46 \\
  concentricCircles1 &    \textbf{53.16} &          54.91 &    55.42 &          54.91 &          54.93 &    55.27 &  \textbf{53.16} &    \textbf{53.16} &    \textbf{53.16} &    59.33 &    62.92 &    55.53 &    55.53 &    53.26 &    56.52 &    54.82 &    53.16 \\
\rowcolor{gray!15}
  concentricCircles2 &   \textbf{153.13} &         154.50 &   156.17 &         161.90 &         161.94 &   163.78 & \textbf{153.13} &            160.00 &            160.00 &   174.40 &   173.13 &   166.06 &   168.97 &   153.59 &   183.92 &   186.95 &   153.13 \\
  concentricCircles3 &   \textbf{270.01} &         274.37 &   275.38 &         274.88 &         274.94 &   275.80 &          271.08 &            323.06 &            324.55 &   297.08 &   295.12 &   306.32 &   304.31 &   271.79 &   336.29 &   387.49 &   270.01 \\
\rowcolor{gray!15}
  concentricCircles4 &   \textbf{451.87} &         459.85 &   461.36 &         479.97 &         480.02 &   487.33 &          454.46 &            510.83 &            522.51 &   520.52 &   527.35 &   513.42 &   511.82 &   455.46 &   502.08 &   686.85 &   451.87 \\
  concentricCircles5 &   \textbf{632.98} &         645.50 &   649.26 &         666.73 &         671.65 &   664.79 &          645.38 &            687.54 &            694.01 &   711.91 &   695.38 &   727.25 &   727.61 &   646.71 &   725.58 &   954.90 &   632.98 \\
\rowcolor{gray!15}
       chaoSingleDep &          1,039.61 &       1,022.98 & 1,023.71 &       1,022.98 &       1,023.05 & 1,022.98 &        1,039.61 & \textbf{1,022.88} & \textbf{1,022.88} & 1,083.60 & 1,084.70 & 1,073.80 & 1,070.00 & 1,042.12 & 1,270.64 & 1,386.75 & 1,022.88 \\
\hline
             Average &            592.20 &         626.94 &   629.29 &         638.14 &         639.28 &   643.11 &          645.10 &            686.86 &            690.91 &   703.51 &   699.35 &   704.96 &   703.12 &   690.58 &   790.44 &   805.27 &   591.03 \\
\hline
                Wins &                44 &              1 &        0 &              1 &              1 &        0 &               7 &                 7 &                 7 &        0 &        0 &        0 &        0 &        0 &        0 &        0 &       48 \\
\hline
         Avg Gap(\%) &              0.13 &           5.01 &     5.56 &           6.20 &           6.37 &     7.12 &           11.40 &             12.36 &             13.08 &    16.56 &    16.60 &    16.73 &    17.04 &    25.77 &    27.85 &    28.59 &     0.00 \\
\hline
Avg Gap(\%) (Losses) &              1.51 &           5.12 &     5.56 &           6.33 &           6.51 &     7.12 &           13.34 &             14.47 &             15.31 &    16.56 &    16.60 &    16.73 &    17.04 &    25.77 &    27.85 &    28.59 &       -- \\
\hline
\end{tabular}}}
{The best solutions are highlighted in bold.}
\end{table}

\begin{table}[ht!]
\TABLE
{Tour lengths obtained by ILS-BCD and by the 10 heuristics evaluated in~\citet{mennell2009thesis} on the 14 benchmark instances with variable radii. \label{tab_cetsp4.1}}
{\resizebox{\fullwidth}{!}{
\begin{tabular}{ccccccccccccc}
\hline
& & \multicolumn{10}{c}{Mennell (2009), Table C1 (p. 232)} & \\
\cmidrule(lr){3-12}
Instance & ILS-BCD & HYBRID2 & HYBRID1 & SZ2 & SZ3$_{50-59}$ & SZ3 & GTSP2 & SZ1 & YUAN+ & LK-SOCP & GTSP1 & Best \\
\hline
\hline
             kroD100 &   \textbf{141.83} &   144.27 &   151.52 &   144.39 &   145.00 &          143.28 & \textbf{141.83} &   145.81 &   149.20 &   149.20 &   142.20 &   141.83 \\
\rowcolor{gray!15}
              rat195 &    \textbf{68.22} &    68.81 &    68.84 &    68.66 &    68.27 &  \textbf{68.22} &           70.18 &    73.51 &    76.79 &    78.89 &   129.71 &    68.22 \\
              lin318 & \textbf{2,078.28} & 2,160.85 & 2,164.45 & 2,175.72 & 2,136.45 &        2,187.68 &        2,081.66 & 2,310.38 & 2,409.20 & 2,365.44 & 2,165.26 & 2,078.28 \\
\rowcolor{gray!15}
               rd400 & \textbf{1,243.39} & 1,252.38 & 1,260.77 & 1,276.08 & 1,270.04 &        1,273.15 &        1,257.01 & 1,279.46 & 1,275.86 & 1,288.05 & 1,258.15 & 1,243.39 \\
              pcb442 &   \textbf{219.84} &   235.19 &   237.32 &   243.48 &   237.83 &          246.08 &          251.22 &   255.97 &   292.72 &   298.99 &   267.03 &   219.84 \\
\rowcolor{gray!15}
                d493 &   \textbf{134.26} &   140.12 &   140.69 &   142.08 &   153.29 &          157.17 &          160.44 &   151.19 &   180.32 &   180.57 &   175.95 &   134.26 \\
             dsj1000 &   \textbf{637.00} &   655.95 &   660.39 &   653.37 &   715.97 &          736.17 &          789.36 &   712.49 &   784.72 &   759.36 &   901.09 &   637.00 \\
\rowcolor{gray!15}
           team1-100 &   \textbf{388.54} &   390.34 &   390.78 &   390.95 &   389.10 & \textbf{388.54} &          388.55 &   410.52 &   413.87 &   413.87 &   392.55 &   388.54 \\
           team2-200 &   \textbf{613.66} &   642.81 &   643.12 &   655.81 &   644.57 &          637.07 &          622.74 &   667.41 &   684.58 &   676.27 &   627.45 &   613.66 \\
\rowcolor{gray!15}
           team3-300 &   \textbf{378.09} &   399.39 &   400.02 &   396.61 &   381.83 &          398.29 &          441.65 &   423.98 &   435.23 &   418.53 &   496.78 &   378.09 \\
           team4-400 & \textbf{1,000.13} & 1,026.32 & 1,028.89 & 1,025.83 & 1,011.77 &        1,020.85 &        1,013.30 & 1,037.44 & 1,054.23 & 1,063.87 & 1,016.00 & 1,000.13 \\
\rowcolor{gray!15}
           team5-499 &   \textbf{446.19} &   456.39 &   457.88 &   455.61 &   476.19 &          454.33 &          528.60 &   496.59 &   517.65 &   557.32 &   631.27 &   446.19 \\
           team6-500 &   \textbf{640.23} &   666.15 &   666.64 &   678.03 &   679.67 &          684.95 &          730.96 &   753.70 &   796.80 &   795.03 &   755.79 &   640.23 \\
\rowcolor{gray!15}
           bonus1000 &   \textbf{980.03} &   996.31 &   996.97 &   992.61 & 1,015.52 &        1,033.07 &        1,049.79 & 1,061.86 & 1,141.71 & 1,149.31 & 1,087.47 &   980.03 \\
\hline
             Average &            640.69 &   659.66 &   662.02 &   664.23 &   666.11 &          673.49 &          680.52 &   698.59 &   729.49 &   728.19 &   717.62 &   640.69 \\
\hline
                Wins &                14 &        0 &        0 &        0 &        0 &               2 &               1 &        0 &        0 &        0 &        0 &       14 \\
\hline
         Avg Gap(\%) &              0.00 &     3.08 &     3.73 &     3.80 &     4.70 &            5.62 &            8.66 &     9.51 &    15.89 &    16.29 &    21.18 &     0.00 \\
\hline
Avg Gap(\%) (Losses) &                -- &     3.08 &     3.73 &     3.80 &     4.70 &            6.56 &            9.32 &     9.51 &    15.89 &    16.29 &    21.18 &       -- \\
\hline
\end{tabular}}}
{The best solutions are highlighted in bold.}
\end{table}

\begin{table}[ht!]
\TABLE
{Tour lengths obtained by ILS-BCD and by four recent heuristics on the 48 benchmark instances of~\citet{mennell2009thesis} with a common radius. \label{tab_recent_common}}
{\resizebox{0.80\fullwidth}{!}{
\begin{tabular}{cccccccc}
\hline
& & \multicolumn{5}{c}{Lei and Hao (2024), Tables 5--6} & \\
\cmidrule(lr){3-7}
Instance & ILS-BCD & SZVNS & (lb/ub)Alg & GA & MA-CETSP & MA-CETSP (avg) & Best \\
\hline
\hline
             kroD100 &   \textbf{159.04} &   \textbf{159.04} &            160.09 &   \textbf{159.04} &   \textbf{159.04} &            159.11 &    159.04 \\
\rowcolor{gray!15}
              rat195 &   \textbf{157.65} &            160.06 &            166.51 &            158.32 &            157.97 &            158.85 &    157.65 \\
              lin318 &          2,830.56 &          2,842.32 &          2,902.53 &          2,838.54 & \textbf{2,816.58} &          2,826.16 &  2,816.58 \\
\rowcolor{gray!15}
               rd400 &          1,026.19 &          1,039.77 &          1,085.75 &          1,032.04 & \textbf{1,018.46} &          1,021.19 &  1,018.46 \\
              pcb442 &            320.66 &            325.02 &            377.51 &            322.54 &   \textbf{319.85} &            321.05 &    319.85 \\
\rowcolor{gray!15}
                d493 &            202.45 &            205.74 &            205.39 &            202.23 &   \textbf{199.02} &            199.67 &    199.02 \\
             dsj1000 &            969.35 &            943.83 &            955.57 &            938.71 &   \textbf{909.50} &            914.21 &    909.50 \\
\rowcolor{gray!15}
             kroD100 &    \textbf{89.67} &             89.92 &    \textbf{89.67} &    \textbf{89.67} &    \textbf{89.67} &             89.68 &     89.67 \\
              rat195 &    \textbf{67.99} &             68.32 &             68.14 &             68.14 &    \textbf{67.99} &    \textbf{67.99} &     67.99 \\
\rowcolor{gray!15}
              lin318 & \textbf{1,394.63} &          1,414.66 &          1,405.07 &          1,404.89 & \textbf{1,394.63} &          1,399.96 &  1,394.63 \\
               rd400 &            466.08 &            474.78 &            460.21 &            458.41 &   \textbf{452.83} &            453.96 &    452.83 \\
\rowcolor{gray!15}
              pcb442 &            142.72 &            152.73 &            146.03 &            145.82 &   \textbf{142.57} &            144.30 &    142.57 \\
                d493 &            100.80 &            102.92 &   \textbf{100.72} &            101.31 &   \textbf{100.72} &            101.12 &    100.72 \\
\rowcolor{gray!15}
             dsj1000 &            377.62 &            393.06 &            374.06 &            376.87 &   \textbf{373.76} &            374.52 &    373.76 \\
             kroD100 &    \textbf{58.54} &    \textbf{58.54} &    \textbf{58.54} &    \textbf{58.54} &    \textbf{58.54} &    \textbf{58.54} &     58.54 \\
\rowcolor{gray!15}
              rat195 &    \textbf{45.70} &    \textbf{45.70} &    \textbf{45.70} &    \textbf{45.70} &    \textbf{45.70} &    \textbf{45.70} &     45.70 \\
              lin318 &            765.98 &            766.16 &   \textbf{765.96} &   \textbf{765.96} &   \textbf{765.96} &   \textbf{765.96} &    765.96 \\
\rowcolor{gray!15}
               rd400 &   \textbf{224.84} &            224.98 &   \textbf{224.84} &   \textbf{224.84} &   \textbf{224.84} &   \textbf{224.84} &    224.84 \\
              pcb442 &    \textbf{83.54} &             83.80 &    \textbf{83.54} &    \textbf{83.54} &    \textbf{83.54} &    \textbf{83.54} &     83.54 \\
\rowcolor{gray!15}
                d493 &    \textbf{69.76} &             69.90 &             69.79 &    \textbf{69.76} &    \textbf{69.76} &    \textbf{69.76} &     69.76 \\
             dsj1000 &            199.97 &            203.07 &   \textbf{199.95} &   \textbf{199.95} &   \textbf{199.95} &   \textbf{199.95} &    199.95 \\
\rowcolor{gray!15}
          team1-100 &   \textbf{307.34} &   \textbf{307.34} &   \textbf{307.34} &   \textbf{307.34} &   \textbf{307.34} &   \textbf{307.34} &    307.34 \\
          team2-200 &   \textbf{246.68} &            246.69 &   \textbf{246.68} &   \textbf{246.68} &   \textbf{246.68} &   \textbf{246.68} &    246.68 \\
\rowcolor{gray!15}
          team3-300 &   \textbf{461.89} &            465.80 &            476.43 &            464.20 &   \textbf{461.89} &            462.16 &    461.89 \\
          team4-400 &            680.09 &            698.05 &            702.69 &            685.52 &   \textbf{669.91} &            675.46 &    669.91 \\
\rowcolor{gray!15}
          team5-499 &            708.85 &            703.38 &            708.45 &            700.50 &   \textbf{693.80} &            697.37 &    693.80 \\
          team6-500 &            225.24 &            226.18 &   \textbf{225.22} &   \textbf{225.22} &   \textbf{225.22} &   \textbf{225.22} &    225.22 \\
\rowcolor{gray!15}
           bonus1000 &            397.19 &            403.06 &            387.13 &            411.85 &   \textbf{384.37} &            390.40 &    384.37 \\
   rotatingDiamonds1 &    \textbf{32.39} &    \textbf{32.39} &    \textbf{32.39} &    \textbf{32.39} &    \textbf{32.39} &    \textbf{32.39} &     32.39 \\
\rowcolor{gray!15}
   rotatingDiamonds2 &   \textbf{140.48} &   \textbf{140.48} &   \textbf{140.48} &   \textbf{140.48} &   \textbf{140.48} &   \textbf{140.48} &    140.48 \\
   rotatingDiamonds3 &   \textbf{380.88} &            380.89 &            380.89 &   \textbf{380.88} &   \textbf{380.88} &   \textbf{380.88} &    380.88 \\
\rowcolor{gray!15}
   rotatingDiamonds4 &   \textbf{770.66} &            770.68 &            772.00 &   \textbf{770.66} &   \textbf{770.66} &            770.67 &    770.66 \\
   rotatingDiamonds5 & \textbf{1,510.75} &          1,510.88 &          1,531.74 & \textbf{1,510.75} & \textbf{1,510.75} &          1,510.77 &  1,510.75 \\
\rowcolor{gray!15}
            bubbles1 &   \textbf{349.13} &            349.14 &            349.14 &            349.14 &            349.14 &            349.14 &    349.13 \\
            bubbles2 &   \textbf{428.28} &   \textbf{428.28} &   \textbf{428.28} &   \textbf{428.28} &   \textbf{428.28} &   \textbf{428.28} &    428.28 \\
\rowcolor{gray!15}
            bubbles3 &   \textbf{529.95} &            532.21 &            529.96 &            529.96 &            529.96 &            529.96 &    529.95 \\
            bubbles4 &            803.39 &            825.33 &            805.56 &            805.46 &   \textbf{802.97} &            805.37 &    802.97 \\
\rowcolor{gray!15}
            bubbles5 &          1,035.82 &          1,073.43 &          1,061.64 &          1,038.16 & \textbf{1,035.32} &          1,037.36 &  1,035.32 \\
            bubbles6 &          1,227.25 &          1,263.68 &          1,313.02 &          1,229.66 & \textbf{1,220.07} &          1,227.15 &  1,220.07 \\
\rowcolor{gray!15}
            bubbles7 &          1,578.31 &          1,639.33 &          1,650.04 &          1,607.31 & \textbf{1,575.04} &          1,586.78 &  1,575.04 \\
            bubbles8 &          1,938.89 &          1,972.99 &          2,021.27 &          1,946.72 & \textbf{1,881.93} &          1,895.74 &  1,881.93 \\
\rowcolor{gray!15}
            bubbles9 &          2,317.46 &          2,330.31 &          2,413.31 &          2,259.22 & \textbf{2,148.40} &          2,168.21 &  2,148.40 \\
  concentricCircles1 &    \textbf{53.16} &    \textbf{53.16} &    \textbf{53.16} &    \textbf{53.16} &    \textbf{53.16} &    \textbf{53.16} &     53.16 \\
\rowcolor{gray!15}
  concentricCircles2 &   \textbf{153.13} &            154.88 &            154.81 &   \textbf{153.13} &   \textbf{153.13} &            153.22 &    153.13 \\
  concentricCircles3 &   \textbf{270.01} &            272.49 &            272.69 &            270.04 &   \textbf{270.01} &            270.04 &    270.01 \\
\rowcolor{gray!15}
  concentricCircles4 &   \textbf{451.87} &            461.36 &            466.52 &            452.64 &   \textbf{451.87} &            453.09 &    451.87 \\
  concentricCircles5 &   \textbf{632.98} &            647.84 &            659.36 &            632.99 &   \textbf{632.98} &            633.43 &    632.98 \\
\rowcolor{gray!15}
       chaoSingleDep & \textbf{1,039.61} &          1,039.63 & \textbf{1,039.61} & \textbf{1,039.61} & \textbf{1,039.61} & \textbf{1,039.61} &  1,039.61 \\
\hline
             Average &            592.20 &            599.05 &            605.74 &            592.02 &            584.31 &            586.47 &    584.31 \\
\hline
                Wins &                27 &                 8 &                16 &                21 &                45 &                17 &        48 \\
\hline
         Avg Gap(\%) &              0.68 &              1.73 &              2.17 &              0.78 &              0.00 &              0.25 &      0.00 \\
\hline
Avg Gap(\%) (Losses) &              1.56 &              2.08 &              3.26 &              1.38 &              0.07 &              0.39 &        -- \\
\hline
        Avg time (s) &          18,469$^a$ &               158 &               460 &             2,762 &               386 &               391 &        -- \\
\hline
\end{tabular}}}
{Values of SZVNS~\citep{Wang2019}, (lb/ub)Alg~\citep{Carrabs2020}, GA~\citep{diplacido2022cetsp}, and MA-CETSP~\citep{Lei2024} (best and average of 20 runs) as reported in~\citet{Lei2024}. Column ``Best'' gives the smallest value among ILS-BCD and the best results of the four methods. Times of the reference methods are those reported in~\citet{Lei2024}, who scaled them to their AMD Opteron 4184 processor using PassMark CPU scores; the number marked with $^a$ in the last row, second column, is the mean CPU time until the best solution was found, on the machine described in Section~\ref{sec4} (not scaled).}
\end{table}

\begin{table}[ht!]
\TABLE
{Tour lengths obtained by ILS-BCD and by four recent heuristics on the 14 benchmark instances of~\citet{mennell2009thesis} with variable radii. \label{tab_recent_var}}
{\resizebox{\fullwidth}{!}{
\begin{tabular}{cccccccc}
\hline
& & \multicolumn{5}{c}{Lei and Hao (2024), Table 7} & \\
\cmidrule(lr){3-7}
Instance & ILS-BCD & SZVNS & (lb/ub)Alg & GA & MA-CETSP & MA-CETSP (avg) & Best \\
\hline
\hline
             kroD100 &   \textbf{141.83} &   \textbf{141.83} &            142.39 &   \textbf{141.83} &   \textbf{141.83} &            141.88 &    141.83 \\
\rowcolor{gray!15}
              rat195 &    \textbf{68.22} &    \textbf{68.22} &    \textbf{68.22} &    \textbf{68.22} &    \textbf{68.22} &    \textbf{68.22} &     68.22 \\
              lin318 &          2,078.28 &          2,082.25 &          2,055.77 &          2,047.42 & \textbf{2,047.11} & \textbf{2,047.11} &  2,047.11 \\
\rowcolor{gray!15}
               rd400 &          1,243.39 &          1,257.73 &          1,305.46 &          1,246.69 & \textbf{1,238.28} &          1,240.89 &  1,238.28 \\
              pcb442 &            219.84 &            221.16 &            220.44 &            220.00 &   \textbf{219.22} &            219.58 &    219.22 \\
\rowcolor{gray!15}
                d493 &            134.26 &            135.02 &            134.74 &            134.28 &   \textbf{134.23} &   \textbf{134.23} &    134.23 \\
             dsj1000 &            637.00 &            625.25 &            625.92 &            624.75 &   \textbf{624.60} &            625.00 &    624.60 \\
\rowcolor{gray!15}
          team1-100 &   \textbf{388.54} &   \textbf{388.54} &   \textbf{388.54} &   \textbf{388.54} &   \textbf{388.54} &   \textbf{388.54} &    388.54 \\
          team2-200 &   \textbf{613.66} &            626.90 &            616.82 &   \textbf{613.66} &   \textbf{613.66} &            614.28 &    613.66 \\
\rowcolor{gray!15}
          team3-300 &   \textbf{378.09} &            379.84 &            378.51 &   \textbf{378.09} &   \textbf{378.09} &   \textbf{378.09} &    378.09 \\
          team4-400 &          1,000.13 &          1,006.71 &          1,025.76 &          1,000.31 &   \textbf{984.24} &            986.44 &    984.24 \\
\rowcolor{gray!15}
          team5-499 &   \textbf{446.19} &   \textbf{446.19} &            446.51 &   \textbf{446.19} &   \textbf{446.19} &   \textbf{446.19} &    446.19 \\
          team6-500 &            640.23 &            621.99 &            626.18 &            620.98 &   \textbf{620.89} &            620.94 &    620.89 \\
\rowcolor{gray!15}
           bonus1000 &            980.03 &            938.27 &            955.41 &            932.07 &   \textbf{917.61} &            921.26 &    917.61 \\
\hline
             Average &            640.69 &            638.56 &            642.19 &            633.07 &            630.19 &            630.90 &    630.19 \\
\hline
                Wins &                 6 &                 4 &                 2 &                 6 &                14 &                 6 &        14 \\
\hline
         Avg Gap(\%) &              1.13 &              0.87 &              1.23 &              0.31 &              0.00 &              0.09 &      0.00 \\
\hline
Avg Gap(\%) (Losses) &              1.97 &              1.22 &              1.44 &              0.54 &                -- &              0.15 &        -- \\
\hline
        Avg time (s) &          24,447$^a$ &                22 &               620 &               330 &               244 &               195 &        -- \\
\hline
\end{tabular}}}
{Values of SZVNS~\citep{Wang2019}, (lb/ub)Alg~\citep{Carrabs2020}, GA~\citep{diplacido2022cetsp}, and MA-CETSP~\citep{Lei2024} (best and average of 20 runs) as reported in~\citet{Lei2024}. Column ``Best'' gives the smallest value among ILS-BCD and the best results of the four methods. Times of the reference methods are those reported in~\citet{Lei2024}, who scaled them to their AMD Opteron 4184 processor using PassMark CPU scores; the number marked with $^a$ in the last row, second column, is the mean CPU time until the best solution was found, on the machine described in Section~\ref{sec4} (not scaled).}
\end{table}

\begin{table}[ht!]
\TABLE
{Best and average tour lengths over 20 runs of the original MA-CETSP (SOCP solved by Gurobi) and of MA-CETSP-BCD, both run by us with the same seeds on the same machine, and best tour length reported in~\citet{Lei2024}. \label{tab_ma_instances}}
{\scriptsize
\setlength{\tabcolsep}{4pt}
\renewcommand{\arraystretch}{0.88}
\resizebox{0.62\fullwidth}{!}{
\begin{tabular}{lcccccc}
\hline
& & \multicolumn{2}{c}{MA-CETSP (Gurobi)} & \multicolumn{2}{c}{MA-CETSP-BCD} & \\
\cmidrule(lr){3-4} \cmidrule(lr){5-6}
Instance & OR & Best & Average & Best & Average & Lei and Hao (2024) \\
\hline
\multicolumn{7}{l}{\textit{Common radius, overlap ratio 0.02}} \\
\hline
kroD100 & 0.02 & \textbf{159.04} & 159.06 & \textbf{159.04} & 159.04 & \textbf{159.04} \\
\rowcolor{gray!15}
rat195 & 0.02 & 157.97 & 158.61 & \textbf{157.87} & 158.38 & 157.97 \\
lin318 & 0.02 & 2,817.20 & 2,829.01 & \textbf{2,815.01} & 2,819.55 & 2,816.58 \\
\rowcolor{gray!15}
rd400 & 0.02 & 1,018.92 & 1,021.36 & \textbf{1,016.45} & 1,019.09 & 1,018.46 \\
pcb442 & 0.02 & 320.11 & 320.47 & \textbf{319.63} & 320.15 & 319.85 \\
\rowcolor{gray!15}
d493 & 0.02 & 199.39 & 199.83 & \textbf{198.15} & 198.78 & 199.02 \\
dsj1000 & 0.02 & 910.71 & 913.60 & \textbf{905.84} & 907.99 & 909.50 \\
\hline
\multicolumn{7}{l}{\textit{Common radius, overlap ratio 0.10}} \\
\hline
kroD100$^*$ & 0.10 & \textbf{89.67} & 89.68 & \textbf{89.67} & 89.67 & \textbf{89.67} \\
\rowcolor{gray!15}
rat195$^*$ & 0.10 & \textbf{67.99} & 68.00 & \textbf{67.99} & 67.99 & \textbf{67.99} \\
lin318$^*$ & 0.10 & \textbf{1,394.63} & 1,403.36 & \textbf{1,394.63} & 1,395.49 & \textbf{1,394.63} \\
\rowcolor{gray!15}
rd400 & 0.10 & 452.68 & 454.94 & \textbf{448.05} & 450.95 & 452.83 \\
pcb442 & 0.10 & 142.70 & 145.25 & \textbf{142.57} & 142.68 & \textbf{142.57} \\
\rowcolor{gray!15}
d493$^*$ & 0.10 & \textbf{100.72} & 100.93 & \textbf{100.72} & 100.72 & \textbf{100.72} \\
dsj1000 & 0.10 & \textbf{373.76} & 374.67 & \textbf{373.76} & 373.76 & \textbf{373.76} \\
\hline
\multicolumn{7}{l}{\textit{Common radius, overlap ratio 0.30}} \\
\hline
kroD100$^*$ & 0.30 & \textbf{58.54} & 58.54 & \textbf{58.54} & 58.54 & \textbf{58.54} \\
\rowcolor{gray!15}
rat195$^*$ & 0.30 & \textbf{45.70} & 45.70 & \textbf{45.70} & 45.70 & \textbf{45.70} \\
lin318$^*$ & 0.30 & \textbf{765.96} & 765.96 & \textbf{765.96} & 765.96 & \textbf{765.96} \\
\rowcolor{gray!15}
rd400$^*$ & 0.30 & \textbf{224.84} & 224.84 & \textbf{224.84} & 224.84 & \textbf{224.84} \\
pcb442$^*$ & 0.30 & \textbf{83.54} & 83.54 & \textbf{83.54} & 83.54 & \textbf{83.54} \\
\rowcolor{gray!15}
d493$^*$ & 0.30 & \textbf{69.76} & 69.76 & \textbf{69.76} & 69.76 & \textbf{69.76} \\
dsj1000$^*$ & 0.30 & \textbf{199.95} & 199.95 & \textbf{199.95} & 199.99 & \textbf{199.95} \\
\hline
\multicolumn{7}{l}{\textit{Common radius, varied overlap ratios}} \\
\hline
team1-100$^*$ & -- & \textbf{307.34} & 307.34 & \textbf{307.34} & 307.34 & \textbf{307.34} \\
\rowcolor{gray!15}
team2-200$^*$ & -- & \textbf{246.68} & 247.55 & \textbf{246.68} & 246.68 & \textbf{246.68} \\
team3-300 & -- & \textbf{461.89} & 462.20 & \textbf{461.89} & 461.89 & \textbf{461.89} \\
\rowcolor{gray!15}
team4-400 & -- & 667.75 & 674.69 & \textbf{663.10} & 666.53 & 669.91 \\
team5-499 & -- & 693.82 & 696.42 & \textbf{693.74} & 695.64 & 693.80 \\
\rowcolor{gray!15}
team6-500$^*$ & -- & \textbf{225.22} & 225.22 & \textbf{225.22} & 225.22 & \textbf{225.22} \\
bonus1000 & -- & 388.31 & 392.96 & \textbf{375.66} & 376.90 & 384.37 \\
\rowcolor{gray!15}
rotatingDiamonds1$^*$ & -- & \textbf{32.39} & 32.39 & \textbf{32.39} & 32.39 & \textbf{32.39} \\
rotatingDiamonds2$^*$ & -- & \textbf{140.48} & 140.48 & \textbf{140.48} & 140.48 & \textbf{140.48} \\
\rowcolor{gray!15}
rotatingDiamonds3 & -- & \textbf{380.88} & 380.88 & \textbf{380.88} & 380.88 & \textbf{380.88} \\
rotatingDiamonds4 & -- & \textbf{770.66} & 770.66 & \textbf{770.66} & 770.66 & \textbf{770.66} \\
\rowcolor{gray!15}
rotatingDiamonds5 & -- & \textbf{1,510.75} & 1,510.78 & \textbf{1,510.75} & 1,510.75 & \textbf{1,510.75} \\
bubbles1$^*$ & -- & \textbf{349.13} & 349.13 & \textbf{349.13} & 349.13 & 349.14 \\
\rowcolor{gray!15}
bubbles2$^*$ & -- & \textbf{428.28} & 428.28 & \textbf{428.28} & 428.28 & \textbf{428.28} \\
bubbles3$^*$ & -- & \textbf{529.95} & 529.96 & \textbf{529.95} & 529.95 & 529.96 \\
\rowcolor{gray!15}
bubbles4 & -- & 802.96 & 805.26 & \textbf{802.83} & 803.50 & 802.97 \\
bubbles5 & -- & 1,035.32 & 1,036.95 & \textbf{1,035.31} & 1,035.31 & 1,035.32 \\
\rowcolor{gray!15}
bubbles6 & -- & 1,220.45 & 1,227.02 & \textbf{1,219.74} & 1,219.88 & 1,220.07 \\
bubbles7 & -- & 1,574.44 & 1,585.42 & \textbf{1,568.74} & 1,570.46 & 1,575.04 \\
\rowcolor{gray!15}
bubbles8 & -- & 1,882.66 & 1,893.64 & \textbf{1,880.70} & 1,880.98 & 1,881.93 \\
bubbles9 & -- & 2,149.84 & 2,167.82 & \textbf{2,145.12} & 2,146.67 & 2,148.40 \\
\rowcolor{gray!15}
concentricCircles1$^*$ & -- & \textbf{53.16} & 53.16 & \textbf{53.16} & 53.16 & \textbf{53.16} \\
concentricCircles2 & -- & \textbf{153.13} & 153.16 & \textbf{153.13} & 153.13 & \textbf{153.13} \\
\rowcolor{gray!15}
concentricCircles3 & -- & \textbf{270.01} & 270.05 & \textbf{270.01} & 270.01 & \textbf{270.01} \\
concentricCircles4 & -- & \textbf{451.87} & 453.06 & \textbf{451.87} & 452.04 & \textbf{451.87} \\
\rowcolor{gray!15}
concentricCircles5 & -- & \textbf{632.98} & 634.18 & \textbf{632.98} & 633.21 & \textbf{632.98} \\
chaoSingleDep & -- & \textbf{1,039.61} & 1,039.61 & \textbf{1,039.61} & 1,039.61 & \textbf{1,039.61} \\
\hline
\multicolumn{7}{l}{\textit{Variable radii}} \\
\hline
kroD100 & -- & \textbf{141.83} & 141.83 & \textbf{141.83} & 141.83 & \textbf{141.83} \\
\rowcolor{gray!15}
rat195$^*$ & -- & \textbf{68.22} & 68.22 & \textbf{68.22} & 68.22 & \textbf{68.22} \\
lin318 & -- & \textbf{2,047.11} & 2,047.11 & \textbf{2,047.11} & 2,047.11 & \textbf{2,047.11} \\
\rowcolor{gray!15}
rd400 & -- & 1,238.02 & 1,239.38 & \textbf{1,238.01} & 1,238.92 & 1,238.28 \\
pcb442 & -- & \textbf{219.22} & 219.75 & \textbf{219.22} & 219.77 & \textbf{219.22} \\
\rowcolor{gray!15}
d493 & -- & \textbf{134.23} & 134.23 & \textbf{134.23} & 134.23 & \textbf{134.23} \\
dsj1000 & -- & \textbf{624.60} & 625.10 & \textbf{624.60} & 624.60 & \textbf{624.60} \\
\rowcolor{gray!15}
team1-100$^*$ & -- & \textbf{388.54} & 388.54 & \textbf{388.54} & 388.54 & \textbf{388.54} \\
team2-200 & -- & \textbf{613.66} & 614.05 & \textbf{613.66} & 613.66 & \textbf{613.66} \\
\rowcolor{gray!15}
team3-300$^*$ & -- & \textbf{378.09} & 378.09 & \textbf{378.09} & 378.09 & \textbf{378.09} \\
team4-400 & -- & 983.63 & 985.87 & \textbf{983.38} & 984.55 & 984.24 \\
\rowcolor{gray!15}
team5-499 & -- & \textbf{446.19} & 446.19 & \textbf{446.19} & 446.19 & \textbf{446.19} \\
team6-500 & -- & \textbf{620.89} & 620.94 & \textbf{620.89} & 620.89 & \textbf{620.89} \\
\rowcolor{gray!15}
bonus1000 & -- & 917.82 & 920.67 & \textbf{916.38} & 917.49 & 917.61 \\
\hline
Average & & 594.77 & 596.54 & 593.99 & 594.47 & 594.67 \\
\hline
\end{tabular}}}
{In boldface, the smallest of the three best values. OR: overlap ratio. Instances marked with $^*$ have known optimal values~\citep{Coutinho2016}.}
\end{table}

\FloatBarrier

\end{appendices}

\end{document}